\documentclass[11pt]{amsart}

\usepackage{latexsym}
\usepackage{amscd}
\usepackage[all]{xy}
\usepackage[mathscr]{euscript}

\normalsize

\RequirePackage{xspace}
\newcommand{\thought}[1]{}
\renewcommand{\thought}[1]{ \textbf{[#1]}}

\usepackage{enumerate}        % Fancy enumerations;  used for Roman numbering.
\newenvironment{roenumerate}{\begin{enumerate}[\upshape (i)]}{\end{enumerate}}
\newcommand\nc {\newcommand}
\newcommand\rnc{\renewcommand}

\newcount\blopone
\newcount\xone
\newcount\xtwo
\newcount\ytwo

\newtheorem{theorem}{Theorem}[section]
\newtheorem{prop}[theorem]{Proposition}
\newtheorem{com}[theorem]{Comment}
\newtheorem{apl}[theorem]{Application}
\newtheorem{exercise}[theorem]{Exercise}
\newtheorem{redu}[theorem]{Reduction}
\newtheorem{refinement}[theorem]{Refinement}
\newtheorem{summary}[theorem]{Summary}
\newtheorem{importnota}[theorem]{Important Notation}
\newtheorem{prblm}[theorem]{Problem}
\newtheorem{notation}[theorem]{Notation}
\newtheorem{explanation}[theorem]{Explanation}
\newtheorem{defin}[theorem]{Definition}
\newtheorem{caution}[theorem]{Caution}
\newtheorem{remark}[theorem]{Remark}
\newtheorem{reminder}[theorem]{Reminder}
\newtheorem{illustration}[theorem]{Illustration}
\newtheorem{observation}[theorem]{Observation}
\newtheorem{elaboration}[theorem]{Elaboration}
\newtheorem{lemma}[theorem]{Lemma}
\newtheorem{construction}[theorem]{Construction}
\newtheorem{discussion}[theorem]{Discussion}
\newtheorem{corollary}[theorem]{Corollary}
\newtheorem{example}[theorem]{Example}
\newtheorem{conclusion}[theorem]{Conclusion}
\newtheorem{sketch}[theorem]{Sketch}
\newtheorem{triviality}[theorem]{Triviality}
\newtheorem{proto}[theorem]{Prototype Quasifibration}
\newtheorem{cauex}[theorem]{Cautionary Example}
\newtheorem{hypo}[theorem]{Hypothesis}
\newtheorem{subth}{ }[theorem]
\newtheorem{case}{Case}[theorem]
\newtheorem{ssubth}{ }[subth]
\newtheorem{facts}[theorem]{Facts}
\newtheorem{history}[theorem]{Historical Survey}
\newtheorem{proofs}[theorem]{Discussion of the Proofs, Old and New}

\nc\tri[1]{\begin{triviality}
\label{#1}}
\nc\fac[1]{\begin{facts}
\label{#1}
\begin{em}}
\nc\app[1]{\begin{apl}
\label{#1}
\begin{em}}
\nc\skt[1]{\begin{sketch}
\label{#1}
\begin{em}}
\nc\hst[1]{\begin{history}
\label{#1}
\begin{em}}
\nc\pfs[1]{\begin{proofs}
\label{#1}
\begin{em}}
\nc\cas[1]{\begin{case}
\label{#1}
\begin{em}}
\nc\rfn[1]{\begin{refinement}
\label{#1}}
\nc\prt[1]{\begin{proto}
\label{#1}}
\nc\lem[1]{\begin{lemma}
\label{#1}}
\nc\pro[1]{\begin{prop}
\label{#1}}
\nc\thm[1]{\begin{theorem}
\label{#1}}
\nc\dis[1]{\begin{discussion}
\label{#1}
\begin{em}}
\nc\cor[1]{\begin{corollary}
\label{#1}}
\nc\dfn[1]{\begin{defin}
\label{#1}}
\nc\sthm[1]{\begin{subth}
\label{#1}}
\nc\exm[1]{\begin{example}
\label{#1}
\begin{em}}
\nc\obs[1]{\begin{observation}
\label{#1}
\begin{em}}
\nc\plm[1]{\begin{prblm}
\label{#1}
\begin{em}}
\nc\rmk[1]{\begin{remark}
\label{#1}
\begin{em}}
\nc\ela[1]{\begin{elaboration}
\label{#1}
\begin{em}}
\nc\rmd[1]{\begin{reminder}
\label{#1}
\begin{em}}
\nc\ntn[1]{\begin{notation}
\label{#1}
\begin{em}}
\nc\exe[1]{\begin{exercise}
\label{#1}
\begin{em}}
\nc\xpl[1]{\begin{explanation}
\label{#1}
\begin{em}}
\nc\smr[1]{\begin{summary}
\label{#1}
\begin{em}}
\nc\cau[1]{\begin{caution}
\label{#1}
\begin{em}}
\nc\hyp[1]{\begin{hypo}
\label{#1}}
\nc\imn[1]{\begin{importnota}
\label{#1}
\begin{em}}
\nc\rdn[1]{\begin{redu}
\label{#1}
\begin{em}}
\nc\cax[1]{\begin{cauex}
\label{#1}
\begin{em}}
\nc\cmt[1]{\begin{com}
\label{#1}
\begin{em}}
\nc\con[1]{\begin{construction}
\label{#1}
\begin{em}}
\nc\ill[1]{\begin{illustration}
\label{#1}
\begin{em}}
\nc\ssthm[1]{\begin{ssubth}
\label{#1}
\begin{em}}
\nc\cnc[1]{\begin{conclusion}
\label{#1}
\begin{em}}

\nc\elem{\end{lemma}}
\nc\erdn{\end{em}\end{redu}}
\nc\erfn{\end{refinement}}
\nc\eprt{\end{proto}}
\nc\ethm{\end{theorem}}
\nc\ecor{\end{corollary}}
\nc\edfn{\end{defin}}
\nc\esthm{\end{subth}}
\nc\epro{\end{prop}}
\nc\etri{\end{triviality}}
\nc\eexm{\end{em}
\end{example}}
\nc\eobs{\end{em}
\end{observation}}
\nc\ecmt{\end{em}
\end{com}}
\nc\efac{\end{em}
\end{facts}}
\nc\eapp{\end{em}
\end{apl}}
\nc\ermk{\end{em}
\end{remark}}
\nc\eela{\end{em}
\end{elaboration}}
\nc\ermd{\end{em}
\end{reminder}}
\nc\eill{\end{em}
\end{illustration}}
\nc\eplm{\end{em}
\end{prblm}}
\nc\ecas{\end{em}
\end{case}}
\nc\eskt{\end{em}
\end{sketch}}
\nc\ecau{\end{em}
\end{caution}}
\nc\ecax{\end{em}
\end{cauex}}
\nc\eimn{\end{em}
\end{importnota}}
\nc\entn{\end{em}
\end{notation}}
\nc\eexe{\end{em}
\end{exercise}}
\nc\expl{\end{em}
\end{explanation}}
\nc\edis{\end{em}
\end{discussion}}
\nc\econ{\end{em}
\end{construction}}
\nc\esmr{\end{em}
\end{summary}}
\nc\ehst{\end{em}
\end{history}}
\nc\epfs{\end{em}
\end{proofs}}
\nc\ehyp{
\end{hypo}}
\nc\ecnc{\end{em}
\end{conclusion}}
\nc\essthm{\end{em}
\end{ssubth}}

\nc\sst{\scriptstyle}
\newcommand{\comment}[1]{}
\newcommand{\ri}{\longrightarrow}

\newcommand{\zz}{{\mathbb Z}}

\newcommand{\nn}{{\mathbb N}}
\newcommand{\K}{{\mathbf K}}
\newcommand{\D}{{\mathbf D}}
\newcommand{\CC}{{\mathbf C}}

\newcommand{\oo}{\otimes}

\nc\op{^{\hbox{\rm\tiny op}}}
\nc\mth{^{\hbox{\rm\tiny th}}}

\nc\script{\mathscr}
\nc\z{\zeta}
\nc\bc{{\mathbb{BC}}}
\nc\ct{{\script T}}
\nc\cf{{\script F}}
\nc\cg{{\script G}}
\nc\ch{{\script H}}
\nc\ck{{\script K}}
\nc\cl{{\script L}}
\nc\cm{{\script M}}
\nc\cn{{\script N}}
\nc\cv{{\script V}}
\nc\ce{{\script E}}
\nc\cs{{\script S}}
\nc\car{{\script R}}
\nc\cd{{\script D}}
\nc\cc{{\script C}}
\nc\ca{{\script A}}
\nc\ci{{\script I}}
\nc\co{{\script O}}
\nc\cu{{\script U}}
\nc\cx{{\script X}}
\nc\cy{{\script Y}}
\nc\cz{{\script Z}}
\nc\bd{\begin{description}}
\nc\ed{\end{description}}
\nc\ctob{{\script C}at\big(\ci^{op},\ca\big)}
\nc\clim{{\ds\mathop{\rm lim}_{\ds\longleftarrow}}\,}
\nc\climi{\clim_{\!i}\,}
\nc\climn{\clim^{\!n}\,}
\nc\colim{{\ds\mathop{\rm colim}_{\ds\la}}}
\nc\colimj{{\ds\mathop{\rm colim}_{\ds\la}}{}_{j\,}}
\nc\oa{\overline{\ca}}
\nc\s{\sigma}
\nc\ta{\tau}
\nc\os{\overline\sigma}
\nc\ot{\overline\tau}
\nc\T{\Sigma}
\nc\Tm{\Sigma^{-1}}
\nc\de[1]{{\mathop{\rm deg(#1)}}}
\nc\Ad[1]{\mathop{\rm Ad}(#1)}
\nc\ad[1]{\mathop{\rm ad}(#1)}
\nc\kth{{\it K}--theory}
\nc\loc[1]{{\text{\rm Loc(#1)}}}
\nc\coloc[1]{{\text{\rm Coloc}(#1)}}

\def\der #1 {D\left(#1\right)}
\nc\prf{\begin{proof}}
\nc\eprf{\end{proof}}
\nc\ds{\displaystyle}
\nc\Tor{\text{\rm Tor}}

\nc\cb{{\script B}}
\nc\ab{{\script A}b}

\nc\be{\begin{roenumerate}}
\nc\ee{\end{roenumerate}}

\nc\cat[1]{{\script C}at\Big({\big\{#1\big\}}\op\,\,,\,\,\ab\Big)}
\nc\csab{{\script C}at\big(\cs^{op},\ab\big)}
\nc\ctab{{\script C}at\Big({\{\ct^\alpha\}}^{op},\ab\Big)}
\nc\csex{{\script E}x\big(\cs^{op},\ab\big)}
\nc\ctex{{\script E}x\Big({\{\ct^\alpha\}}^{op},\ab\Big)}
\nc\sub{\qquad\subset\qquad}
\nc\ctr[1]{{\left.\ct\left(-,#1\right)\right|}_{\cs}}
\nc\ctrf[2]{{\left.\ct\left(#1,#2\right)\right|}_{\cs}}
\nc\Ctr[1]{{\left.\ct\left(-,#1\right)\right|}_{\ct^\alpha}}
\nc\Ctrf[2]{{\left.\ct\left(#1,#2\right)\right|}_{\ct^\alpha}}

\nc\la{\longrightarrow}
\nc\nin{\noindent}
\nc\cad[1]{\text{card}(#1)}
\nc\eq{\quad=\quad}
\nc\BA{\begin{array}{c}}
\nc\EA{\end{array}}
\nc\barr{
\[
\begin{array}{cccccccccccccccc}
}
\nc\earr{
\end{array}
\]
}
\nc\as[1]{{\langle S\rangle}^{#1}}
\nc\sh{\text{\it shift}}

\nc\yy[1]{{\left.\ct\left(-,#1\right)\right|}_{\ct^c}}
\nc\vrep[2]{{\left.\ct\left(#1,#2\right)\right|}_{\ct^\alpha}}
\nc\da{\downarrow}
\nc\Hom{{\mathop{\rm Hom}}}
\nc\HHom{{\script H}{\mathop{\rm om}}}
\nc\End{{\mathop{\rm End}}}
\nc\Ext{{\mathop{\rm Ext}}}
\nc\mr\Modtc
\nc\PExt{{\mathop{\rm PExt}}}
\nc\stm{\text{\rm stmod}(kG)}
\nc\stM{\text{\rm StMod}(kG)}
\nc\e{\varepsilon}
\nc\p{\varphi}

\nc\rs{\s^{-1}A}
\nc\br{{\{\s^{-1}A\}}}
\nc\ra\ri
\nc\y[1]{\mathbf{y}#1}
\nc\x[1]{\mathbf{z}#1}
\nc\mmod[1]{#1\text{--\rm mod}}
\nc\Mod[1]{#1\text{--\rm Mod}}
\nc\Md {\ensuremath{\mathop{\textup{Mod}}}}
\rnc\mod[1]{\ensuremath{\mathop{#1\textup{--mod}}}\xspace}
\nc\MMod[1]{\text{Mod-}#1}
\nc\Modtc{\Mod{\ct^c}}
\nc\pgldim[1]{\mathop{\rm pgldim}\,#1}
\nc\tf{{\rm [TR5]}}
\nc\tfs{{\rm [TR5$^*$]}}
\nc\Fun{\text{\rm Funct}(F\op,\ab)}
\nc\sym{\text{\rm Sym}}
\nc\sgn{\text{\rm sgn}}
\nc\Pro{\text{\rm Prod}^{}_\alpha(F\op,\ab)}
\nc\Yt[1]{{\left.\Hom_\ct^{}\left(-,#1\right)\right|}_F^{}}
\nc\dl{\delta}
\nc\Proj[1]{#1\text{--\rm Proj}}
\nc\proj[1]{#1\text{--\rm proj}}
\nc\Flat[1]{#1\text{--\rm Flat}}
\nc\Inj[1]{#1\text{--\rm Inj}}
\nc\Ima{\mathrm{Im}}
\nc\Ker{\mathrm{Ker}}
\nc\ov{\overline}
\nc\wi{\wt{\text{\it\i}}}
\nc\wt{\widetilde}
\nc\ph{\varphi}
\nc\tstr{{\it t}--structure}
\nc\spec[1]{{\text{\rm Spec}(#1)}}

\nc\EProd{\text{\rm EProd}}
\nc\ECoprod{\text{\rm ECoprod}}
\nc\Prod{\text{\rm Prod}}
\nc\ldimp{\text{\rm LDim}^{\prod}}
\nc\ldimc{\text{\rm LDim}^{\coprod}}
\nc\gen[2]{{\langle#1\rangle}^{}_{#2}}
\nc\genu[3]{{\langle#1\rangle}^{[#3]}_{#2}}
\nc\ogen[1]{\ov{\langle#1\rangle}}
\nc\ogenun[2]{\ov{\langle#1\rangle}_{#2}^{}}
\nc\ogenu[3]{\ov{\langle#1\rangle}^{[#3]}_{#2}}
\nc\ogenul[3]{\ov{\langle#1\rangle}^{(-\infty,#3]}_{#2}}
\nc\genuf[3]{{\langle#1\rangle}^{[#3,\infty)}_{#2}}
\nc\genul[3]{{\langle#1\rangle}^{(-\infty,#3]}_{#2}}
\nc\dperf[1]{\D^{\mathrm{perf}}(#1)}
\nc\dcoh{\mathbf{D}^b_{\mathrm{coh}}}
\newcommand{\Dqc}{{\mathbf D_{\text{\bf qc}}}}
\newcommand{\Dqcmi}{{\mathbf D_{\text{\bf qc}}^-}}
\newcommand{\Dqcpl}{{\mathbf D_{\text{\bf qc}}^+}}
\newcommand{\Dqcb}{{\mathbf D_{\text{\bf qc}}^b}}

\nc\dmcoh{\mathbf{D}^-_{\mathrm{coh}}}
\nc\dscoh{\mathbf{D}^{}_{\mathrm{coh}}}
\nc\RHHom{{\script{RH}}{\mathrm{om}}}
\nc\Coprod{\mathrm{Coprod}}
\nc\COprod{\mathrm{coprod}}
\nc\add{\mathrm{add}}
\nc\Add{\mathrm{Add}}
\nc\Smr{\mathrm{smd}}
\nc\id{\mathrm{id}}
\nc\LL{\mathbf{L}}
\nc\R{\mathbf{R}}
\nc\fc{\mathfrak{C}}
\nc\fl{\mathfrak{L}}
\nc\fs{\mathfrak{S}}

\nc\hoco{
\begin{picture}(40,10)
\put(20,0){\makebox(0,0)[b]{\text{\rm Hocolim}}}
\put(5,-2){\vector(1,0){30}}
\end{picture}\,\,}

\nc\holim{
\begin{picture}(40,10)
\put(20,0){\makebox(0,0)[b]{\text{\rm Holim}}}
\put(35,-2){\vector(-1,0){30}}
\end{picture}}

\begin{document}

\author{Amnon Neeman}\thanks{The research was partly supported 
by the Australian Research Council}
\address{Centre for Mathematics and its Applications \\
        Mathematical Sciences Institute\\
        Building 145\\
        The Australian National University\\
        Canberra, ACT 2601\\
        AUSTRALIA}
\email{Amnon.Neeman@anu.edu.au}

\title{Approximable triangulated categories}

\begin{abstract}
In this survey we present the relatively
new concept of \emph{approximable triangulated categories.}
We will show that the definition is natural, that it leads to powerful
new results, and that it throws new light on old, familiar objects.
\end{abstract}



\maketitle

\tableofcontents

\section{Introduction}
\label{S0}

In this survey there is one major, recent technical tool---we present
the concept of \emph{approximable triangulated categories.}
And then we will sketch some results, from the last few years, showing
that this tool is useful.

The definition of approximable triangulated categories relies on
the following building blocks: compact
generators in triangulated categories and {\it t}--structures.
We have an extensive background section introducing
these---we recommend that beginners skip the remainder of the 
introduction and move on to Section~\ref{S1}.
The introduction of an article is normally the author's attempt to
persuade the expert to read on---hence it tends to assume some
familiarity with the existing theory, the expert will not want to be
bored with stuff she already knows, she will want to see if this
article contains anything new and interesting. The introduction to
an article is often more intimidating than the body of the manuscript.

In addition to Section~\ref{S1}
the beginners might wish to look at the three appendices, which
were written to answer questions from students to whom the material
was new. Especially relevant is Appendix~\ref{A3},
where we draw up a dictionary between our approach and the more
standard one in the literature---in the interest of efficiency we
depart from the usual way to introduce derived categories. An
interested
beginner, who wants to explore this further by
looking elsewhere, is encouraged to consult this dictionary.

Back to the experts: we plan to discuss approximability
in triangulated categories and its applications, and we begin
with a heuristic explanation of what approximability is all
about.

Any \tstr\ on the triangulated category can be used to define a ``metric'':
two objects are close to each other if they agree up to a
small ``difference''. 
More precisely: the objects $x,y\in\ct$ are close to each other if
there exists in $\ct$ a triangle $x\la y\la z\la$ with
$z\in\ct^{\leq-n}$ for some large $n$.\footnote{The category $\ct$ 
isn't quite a metric space,
the obvious ``metric'' isn't symmetric. There is a rich literature on
asymmetric metrics
following Lawvere~\cite{Lawvere73,Lawvere02}, but the emphasis 
there 
goes in entirely different directions---it's the triangle inequality
in
Lawvere's asymmetric metrics that is emphasized, see for example 
Leinster~\cite{Leinster13}. Our metrics might be asymmetric, but they
are decidedly non-archimedean---hence the triangle inequality is dumb,
every triangle is isosceles.}
We declare that, the larger the integer $n$, the closer
the points $x$ and $y$.
In the world of metric spaces we are accustomed
to the notion of equivalent metrics,
and this naturally leads to the concept of equivalent {\it t}--structures.

We are also accustomed to expressing points in a metric
space as limits of
Cauchy sequences of simpler, more accessible points. For
example the Taylor series approximates a function by polynomials,
and the Fourier series approximates a function by finite sums
of exponentials. There is
a triangulated category version, we will explain it more fully
in the body of the article.
For the introduction the discussion below will give the gist of
the construction, albeit a little vaguely and with details missing.

\dis{D0.1}
If we plan to approximate objects of the triangulated
category $\ct$, by Cauchy sequences of simpler objects,
then we first need to measure
what we mean by ``simplicity''---returning to the analogy
of the previous
paragraphs, we need to declare what will be
the triangulated category replacement for the polynomials
forming the partial sums in a
Taylor series. In doing this we
will slightly modify an idea due to Bondal and
Van~den~Bergh~\cite{BondalvandenBergh04}. We will start with a
compact generator $G$ for the triangulated category and, for
each integer $n>0$, we will define two classes of objects
$\genu Gn{-n,n}\subset\ogenu Gn{-n,n}$. These will be the objects
obtainable from $G$ in $n$ allowable steps---the difference between
$\genu Gn{-n,n}$ and $\ogenu Gn{-n,n}$ is that for
the smaller $\genu Gn{-n,n}$ there are fewer operations allowed.

Returning to the analogy with Taylor series:
so far we have explained what will be our replacement for
the polynomials of degree $\leq n$. We have already indicated the ``metric''
we plan to work with,
it is the one determined by whatever \tstr\ we end up choosing.
So it becomes interesting to know which objects in the
triangulated category have Taylor series ``converging'' to them. And now
we come to the (somewhat imprecise)
definition: the triangulated
category $\ct$ is declared to be \emph{approximable} if it has coproducts,
and there exist in $\ct$
\be
\item
a compact generator $G$
\item
a \tstr\ $(\ct^{\leq0},\ct^{\geq0})$
\setcounter{enumiv}{\value{enumi}}
\ee
and these \tstr\ and generator can be chosen to satisfy
\be
\setcounter{enumi}{\value{enumiv}}
\item
For some $n>0$ we have $G[n]\in\ct^{\leq0}$ and
$\Hom\big(G[-n]\,\,,\,\,\ct^{\leq0}\big)=0$.
\item
In the metric induced by the \tstr\ $(\ct^{\leq0},\ct^{\geq0})$ of (ii), every
object in $\ct^{\leq0}$ can be expressed as the limit of a
sequence whose terms belong to $\cup_n\ogenu Gn{-n,n}$.
\ee
\edis

\rmk{R0.3}
It is a formal consequence of the definition that an approximable
triangulated category is complete with respect to the metric
of Discussion~\ref{D0.1}(iv)---any Cauchy
sequence converges. Also:
by definition, if $\ct$ is approximable then $\ct^{\leq0}$ is
contained in the closure of $\cup_n\ogenu Gn{-n,n}$ with
respect to the metric. It turns out that the closure is not much larger:
it is nothing other than $\ct^-=\cup_n\ct^{\leq n}$.

One could also wonder what the closure of $\cup_n\genu Gn{-n,n}$
might be---we will return to this later, it turns out to be
a subcategory which, for many $\ct$, has been studied extensively
in the classical literature.
\ermk

Now that we have a rough idea what approximability means, the first
question we might ask ourselves is ``Is the theory nonempty, are
there any examples?''

\exm{E0.5}
It turns out there are plenty of examples. If $R$ is any ring
then $\D(R)$, the unbounded derived category of complexes of
left $R$--modules, is
an example. So is the homotopy category of spectra, and so is
the category $\D(\Mod R)$, provided $R$ is a dga such that
$H^i(R)=0$ for all $i>0$. Here $\D(\Mod R)$ stands for the derived category
whose objects are all left dg $R$--modules.

All of these are easy examples. It is a nontrivial
theorem that, when $X$ is a quasicompact, separated scheme, the category
$\Dqc(X)$ is approximable. Here $\Dqc(X)$ means the
derived category, whose
objects are cochain complexes of $\co_X^{}$--modules with quasicoherent
cohomology. It is also a nontrivial theorem that, under
reasonable hypotheses, the recollement of two approximable
categories is approximable.
\eexm

\app{A0.7}
We have introduced a new gadget---namely approximable
triangulated categories---and mentioned that there are plenty of
interesting examples out there.
But the skeptical reader will want to know what
use this new toy might have: are there applications? Do we learn anything new,
about the familiar old categories
of Example~\ref{E0.5}, because we now know them to be approximable?

The answer is Yes. We list below
five results we were recently able to prove, using
the fact that $\Dqc(X)$ is approximable.
\be
\item
Suppose $X$ is a quasicompact, separated scheme.
Then the category $\dperf X$ is strongly generated, in the sense
of Bondal and Van~den~Bergh~\cite{BondalvandenBergh04},
if and only if $X$ can be covered by open affine subsets $U_i=\spec{R_i}$
with each $R_i$ of finite global dimension.
\item
Suppose $X$ is a finite-dimensional, noetherian, separated scheme,
and assume further that every closed, reduced, irreducible
subscheme of $X$ has a regular alteration.
Then the category
$\dcoh(X)$ is strongly generated.
\item
Suppose $X$ is a scheme proper over a noetherian ring $R$,
and let $\cy:\dcoh(X)\la\Hom_R^{}\big[\dperf X\op\,\,,\,\,\Mod R\big]$
be the Yoneda map. That is: $\cy$ is the map taking an object
$B\in\dcoh(X)$ to the functor $\cy(B)=\Hom(-,B)$,
viewed as an $R$--linear homological
functor $\dperf X\op\la\Mod R$.

Then $\cy$ is fully faithful, and the essential image
of $\cy$ are the finite homological functors. A functor
$H:\dperf X\op\la \Mod R$ is \emph{finite} if, 
for any object $A\in\dperf X$, the $R$--modules $H\big(A[i]\big)$ are
all finite and all but finitely many of them vanish.
\item
Suppose $X$ is a finite-dimensional scheme proper over a noetherian ring $R$,
and assume further that every closed, reduced, irreducible
subscheme of $X$ has a regular alteration.
Let $\wt\cy:{\dperf X}\op\la\Hom_R^{}\big[\dcoh(X)\,\,,\,\,\Mod R\big]$
be the Yoneda map. That is: $\wt\cy$ is the map taking an object
$A\in\dperf X$ to the functor $\wt\cy(A)=\Hom(A,-)$,
viewed as an $R$--linear homological
functor $\dcoh(X)\la\Mod R$.

Then $\wt\cy$ is fully faithful, and the essential image
of $\wt\cy$ are the finite homological functors.
\item
Suppose $X$ is a noetherian, separated scheme. Then the categories
$\dperf X$ and $\dcoh(X)$ determine each other. More explicitly: there
is a recipe which takes a triangulated category $\cs$ as input,
and out of it cooks up another triangulated category denoted
$\fs(\cs)$. If we apply this recipe to $\dperf X$ what comes out is
$\dcoh(X)$, and if we apply it to $\big[\dcoh(X)\big]\op$ the output
is $\big[\dperf X\big]\op$. 
\setcounter{enumiv}{\value{enumi}}
\ee
In the body of the article we will say more about the theorems---for example
we will remind the reader what it means for a triangulated category $\ct$
to be ``strongly generated''. 
For now we note only that (i), (ii), (iii), (iv) and (v) above represent
sharp improvements over the existing literature. More precisely 
\be
\setcounter{enumi}{\value{enumiv}}
\item
There were versions of (i) and (ii) known before approximability, but
they all assumed equal
characteristic---the reader can find a sample of the
known results 
in Bondal and Van~den~Bergh~\cite[Theorem~3.1.4]{BondalvandenBergh04},
Iyengar and Takahashi~\cite[Corollary~7.2]{Iyengar-Takahashi15},
Orlov~\cite[Theorem~3.27]{Orlov16} and
Rouquier~\cite[Theorem~7.38]{Rouquier08}.
\item
The only known versions of (iii) and (iv), prior to approximability, assumed that $R$ is
a field. See  
Bondal and Van~den~Bergh~\cite[Theorem~A.1]{BondalvandenBergh04} for
(iii), and Rouquier~\cite[Corollary~7.51(ii)]{Rouquier08} for (iv).
\item
  The only known versions of (v) prior to approximability assumed either that
  $X$ is affine, see Rickard~\cite[Theorem~6.4]{Rickard89b}, or that 
$X$ is projective over a field, see 
Rouquier~\cite[Remark~7.50]{Rouquier08}.
\ee
\eapp

The definition of approximability is the assumption that
there exist a \tstr\ $(\ct^{\leq0},\ct^{\geq0})$
and a compact
generator $G\in\ct$ with some properties. It becomes natural to wonder
how free we are in the choice of \tstr\ and compact generator. This leads
to a string of surprising results.

\fac{F0.9}
Let $\ct$ be a triangulated category with coproducts, and assume it has
a compact generator $G$. Then the following can be proved.
\be
\item
There exists a preferred equivalence class of {\it t}--structures
in $\ct$. Here two {\it t}--structures are declared equivalent if they
induce equivalent metrics.
\item
Let us choose in $\ct$ a compact generator $G$ and a \tstr\ 
$(\ct^{\leq0},\ct^{\geq0})$, and assume they satisfy the conditions in
Discussion~\ref{D0.1}~(iii) and (iv)---that is:
the pair $G$ and $(\ct^{\leq0},\ct^{\geq0})$ pass the test
for checking the approximability of $\ct$.

Then it's automatic that
$(\ct^{\leq0},\ct^{\geq0})$ belongs to the preferred equivalence
class of {\it t}--structures. Hence the metric defined
by any \tstr\ $(\ct^{\leq0},\ct^{\geq0})$, which can be used to test
for approximability, is unique up to equivalence.
\item
Suppose $\ct$ is approximable.
Then any \tstr\ $(\ct^{\leq0},\ct^{\geq0})$ in the preferred
equivalence class and any compact generator $G$
satisfy the conditions in
Discussion~\ref{D0.1}~(iii) and (iv).
\setcounter{enumiv}{\value{enumi}}
\ee
Thus approximability is robust; it doesn't really matter which \tstr\
and compact generator one chooses, as long as the \tstr\ belongs
to the preferred equivalence class. Furthermore the
categories
\[
\ct^-=\cup_n\ct^{\leq n}\,,\qquad \ct^+=\cup_n\ct^{\geq -n}\,,
\qquad \ct^b=\ct^-\cap\ct^+
\]
turn out to be intrinsic. They depend only on $\ct$,
not on the particular representative
$(\ct^{\leq0},\ct^{\geq0})$ in the preferred equivalence class.

Now that we know the metric is intrinsic (up to equivalence), it
makes sense to return to the question raised in Remark~\ref{R0.3}.
What is the closure of $\cup_n\genu Gn{-n,n}$? In view of the above
we should not be surprised to learn
\be
\setcounter{enumi}{\value{enumiv}}
\item
Define the category $\ct^-_c$ to be the closure in $\ct$ of
$\cup_n\genu Gn{-n,n}$. This category is intrinsic, it does not depend on
the choice of compact generator $G$. And it follows that
the category $\ct^b_c=\ct^-_c\cap\ct^b$ must also be intrinsic.
\ee
\efac

\rmk{R0.11}
It becomes interesting to figure out what these intrinsic subcategories
are in the special cases of Example~\ref{E0.5}. Let us confine ourselves
to just one case: assume $\ct=\Dqc(X)$ with $X$ a separated, noetherian scheme.
In this special case one can prove:
\be
\item
The standard \tstr\ is in the preferred equivalence class. Hence the categories
$\ct^-$, $\ct^+$ and $\ct^b$ have their usual meanings: that is
$\ct^-=\Dqcmi(X)$, $\ct^+=\Dqcpl(X)$ and $\ct^b=\Dqcb(X)$.
\item
The category $\ct^-_c$ turns out to be $\dmcoh(X)$, hence the category
$\ct^b_c=\ct^-_c\cap\ct^b$ is nothing other than $\dcoh(X)$.
\setcounter{enumiv}{\value{enumi}}
\ee
It turns out that Applications~\ref{A0.7}~(iii), (iv) and (v)
generalize greatly. The glorious,
abstract versions for (iii) and (iv) go as follows.
Let $R$ be a noetherian ring, and let
$\ct$ be an $R$--linear, approximable triangulated category.
Suppose there exists in $\ct$ a compact generator $G$, such that
$\Hom\big(G,G[n]\big)$ is a finite $R$--module for all $n\in\zz$.
Consider the two functors
\[
\cy:\ct^-_c\la\Hom_R^{}\big([\ct^c]\op\,,\,\Mod R\big),\qquad
\wt\cy:\big[\ct^-_c\big]\op\la\Hom_R^{}\big(\ct^b_c\,,\,\Mod R\big)
\]
defined by the formulas $\cy(B)=\Hom(-,B)$ and
$\wt\cy(A)=\Hom(A,-)$. Note that, in these formulas,
we permit all $A,B\in\ct^-_c$. But the $(-)$ in the formula
$\cy(B)=\Hom(-,B)$ is assumed to belong to $\ct^c$,
whereas the $(-)$ in the formula $\wt\cy(A)=\Hom(A,-)$ must
lie in $\ct^b_c$.
Now consider the following composites
\[\xymatrix@C+20pt@R-20pt{
\ct^b_c \ar@{^{(}->}[r]^i & \ct^-_c
\ar[r]^-{\cy} &
\Hom_R^{}\big([\ct^c]\op\,,\,\Mod R\big) \\
\big[\ct^c\big]\op \ar@{^{(}->}[r]^{\wi} & \big[\ct^-_c\big]\op
\ar[r]^-{\wt\cy} &
\Hom_R^{}\big(\ct^b_c\,,\,\Mod R\big)
}\]
where $i$ and $\wi$ are the obvious inclusions.
We assert:
\be
\setcounter{enumi}{\value{enumiv}}
\item
  The functor $\cy$ is full, and the essential image consists
  of the
  locally finite homological functors. A functor
  $H:\big[\ct^c\big]\op\la \Mod R$ is \emph{locally finite} if
$H\big(A[n]\big)$ is a finite $R$--module for every $n\in\zz$, and
vanishes if $n\ll0$. 

The composite $\cy\circ i$
  is fully faithful, and the essential image consists of the
  finite homological functors.
\item
  Assume there exists an integer $N>0$ and an object $G'\in\ct^b_c$
  with $\ct=\ogen {G'}_N^{(-\infty,\infty)}$---this condition will
  be explained in the body of the paper, under the hypotheses
  placed on $X$ in Application~\ref{A0.7}(iv) the condition is
  satisfied by $\ct=\Dqc(X)$,
  this may be found in \cite[Theorem~2.3]{Neeman17}.

  Then
  the functor $\wt\cy$ is full, and the essential image consists
  of the
  locally finite homological functors. The composite $\wt\cy\circ \wi$
  is fully faithful, and the essential image consists of the
  finite homological functors.
\setcounter{enumiv}{\value{enumi}}
\ee
As we have said, Application~\ref{A0.7}(v) also has a vast
generalization, which goes as follows:
\be
\setcounter{enumi}{\value{enumiv}}
\item
With the notation as in Application~\ref{A0.7}(v) one has, for
any approximable $\ct$, a triangulated equivalence
$\fs(\ct^c)\cong\ct^b_c$. If the triangulated category $\ct$ is not only
approximable but also \emph{noetherian,} then one also has a triangulated
equivalence
$\fs\big(\big[\ct^b_c\big]\op\big)\cong
\big[\ct^c\big]\op$.
\ee

The notion of noetherian triangulated categories
in Remark~\ref{R0.11}(v) is new, and was
inspired by the result. Noetherianness is a condition that seems natural,
and guarantees that there will be plenty of objects in
$\ct^b_c$. Without some noetherian hypothesis, the only obvious 
object in $\ct^b_c$ is zero.
\ermk

\medskip

\nin
    {\bf Acknowledgements.}\ \
    The author would like to thank Jesse Burke, Anthony Kling, Steve Lack, 
    Luke Mitchelson, Bregje Pauwels, Geordie Williamson
and an anonymous referee
for questions and comments that led to improvements on earlier versions of 
the manuscript. These comments were made both about earlier
drafts, and during talks
presenting parts or all of the material.

\section{Background}
\label{S1}

It's time to speak to the non-experts---the readers familiar with
triangulated categories, compact generators and {\it t}--structures
are advised to skip ahead to Section~\ref{S2}. In this section we will
present a
quick reminder of the three concepts in the sentence above.
We plan to usually proceed from the concrete to the abstract: for most of the
section we study first an example, actually four examples---all four
examples
will be derived categories  $\D_{\fc}^{\fc'}(\ca)$,
we list the four 
in Examaple~\ref{E1.100372}---and
only then do we move on to the general definitions. We should therefore begin
by recalling what are the categories $\D_{\fc}^{\fc'}(\ca)$,
first in generality that covers the four examples and more, and then
narrowing down to the specific ones that will interest us. 
\exm{E1.1}
Let $\ca$ be an abelian category. The derived category $\D_{\fc}^{\fc'}(\ca)$ is as follows:  
\be 
\item
The objects are the cochain complexes in $\ca$, that is 
diagrams in $\ca$ of the form
\[\xymatrix@C+15pt{
\cdots \ar[r]& A^{-2}  \ar[r]& A^{-1}  \ar[r]& A^{0}  \ar[r]
&  A^{1} \ar[r]& A^{2} \ar[r] &\cdots
}\]
where the composites $A^i\la A^{i+1}\la A^{i+2}$ all vanish. The subscript
$\fc$ and
superscript $\fc'$ stand for conditions.
We may choose not to
allow all cochain complexes, when the mood strikes us we
can capriciously impose any conditions on the cochain complexes
that our heart desires---subject to the mild hypotheses that
guarantee that the few operations we're about to perform take complexes
satisfying the restrictions to complexes satisfying the restrictions.
\item
Cochain maps are morphisms in $\D_{\fc}^{\fc'}(\ca)$, that is any commutative
diagram
\[\xymatrix@C+10pt{
\cdots \ar[r]& A^{-2}  \ar[d]\ar[r]& A^{-1} \ar[d] \ar[r]& A^{0}\ar[d]  \ar[r]
&  A^{1} \ar[d]\ar[r]& A^{2}\ar[d] \ar[r] &\cdots\\
\cdots \ar[r]& B^{-2}  \ar[r]& B^{-1}  \ar[r]& B^{0}  \ar[r]
&  B^{1} \ar[r]& B^{2} \ar[r] &\cdots
}\]
is a morphism from the top to the bottom
row---as long as the rows are cochain complexes satisfying the
restrictions, that is objects in $\D_{\fc}^{\fc'}(\ca)$.

But then we formally invert the cohomology isomorphisms. In the literature
the cohomology isomorphisms often go by the name ``quasi-isomorphisms''.
\ee 
\eexm

\xpl{X1.3}
Given a category $\cc$ and a collection $S$ of morphisms in $\cc$,
an old theorem of Gabriel and Zisman~\cite{Gabriel-Zisman67}
tells us that there exists
a functor $F:\cc\la S^{-1}\cc$ so that
\be
\item
If $s\in S\subset\text{Mor}(\cc)$ then $F(s)$ is invertible.
\item
Any functor $F'':\cc\la\cb$, with $F''(S)$ contained in the isomorphisms
of $\cb$, factors uniquely as
$\cc\stackrel{F}\la S^{-1}\cc\stackrel{F'}\la \cb$.
\ee
What we mean when we say that in $\D_{\fc}^{\fc'}(\ca)$ we ``formally invert''
the cohomology isomorphisms is: let $\cc$ be the category with the
same objects as
$\D_{\fc}^{\fc'}(\ca)$ but where the morphisms are the cochain maps,
and let $S$ be the collection of cochain maps inducing cohomology isomorphisms.
Then $\D_{\fc}^{\fc'}(\ca)$ is defined to be $S^{-1}\cc$.

In principle categories of the form $S^{-1}\cc$ can be dreadful---the morphisms
are equivalence classes of composable strings, 
where each string is a sequence whose pieces are either morphisms in $\cc$
or inverses of elements of $S$.
The Hom-sets needn't be small, and in general
it can be a nightmare to decide when two such strings are equivalent, meaning
define the same morphism in $S^{-1}\cc$. For
categories like $\D_{\fc}^{\fc'}(\ca)$ the calculus of fractions happens not
to be too bad, there is
a literature dealing with it. The interested reader is referred to
Hartshorne~\cite{Hartshorne66} or Verdier~\cite{Verdier96} for the original
presentations, or Gelfand and Manin~\cite{GelfandManin2},
Kashiwara and Schapira~\cite{Kashiwara-Schapira06} or
Weibel~\cite{Weibel94} for more modern treatments.

In this survey we skip the discussion of the calculus of fractions.
This means that the reader will be asked to believe several
computations along the way---when these occur there will be a footnote
to the effect.
\expl

\exm{E1.100372}
In this survey, the key examples to keep in mind are:
\be
\item
If $R$ is a ring, $\D(R)$ will be our shorthand for $\D(\Mod R)$;
the abelian category $\ca$ is the category of all left
$R$--modules,  
and since there are no superscripts or subscripts decorating 
the $\D$ we impose no conditions.
All cochain complexes of left $R$-modules are objects of $\D(R)$.
\setcounter{enumiv}{\value{enumi}}
\ee
Now let $X$ be a scheme. The abelian category,
in all three examples below, is the category
of sheaves of $\co_X^{}$--modules. It's customary
to abbreviate what should be written $\D(\Mod{\co_X^{}})$
to just $\D(X)$, and we will follow this custom.

But all three categories we will look at are decorated, there are
restrictions. We list them below.
\be
\setcounter{enumi}{\value{enumiv}}
\item
The objects in $\Dqc(X)$ are cochain
complexes of $\co_X^{}$--modules,
and the only condition we impose is that 
the cohomology sheaves must be quasicoherent.
\item
The objects of $\dperf X$ are the perfect complexes.
A cochain complex of $\co_X^{}$--modules is
\emph{perfect} if it is locally isomorphic to a bounded complex
of vector bundles. More precisely: an object
$P\in\Dqc(X)$ belongs to $\dperf X$ if
there exists an open cover of $X$ of
the form $X=\cup_i^{}U_i$ such that, if $u_i:U_i\la X$ is the
inclusion, then the obvious functor $u^*_i:\Dqc(X)\la\Dqc(U_i)$
takes $P\in\Dqc(X)$ to an object $u_i^*(P)\in\Dqc(U_i)$ which is isomorphic
in $\Dqc(U_i)$ to a bounded complex of vector bundles.
\item
Assume $X$ is noetherian. The objects of $\dcoh(X)$ are the complexes
of $\co_X^{}$--modules
with coherent cohomology---as indicated by the subscript---and this cohomology
vanishes in all but finitely many degrees, the superscipt 
$b$ stands for ``bounded''.
\ee
\eexm

\rmk{R1.4}
If we're going to be working with categories like $\D_{\fc}^{\fc'}(\ca)$,
it is natural to wonder what useful structure they might have.
The next definition spells out the answer.

The idea is simple enough: we started with the category $\cc$
whose objects are the same as those of $\D_{\fc}^{\fc'}(\ca)$,
but the morphisms were honest cochain maps.
We then performed the construction of Explanation~\ref{X1.3},
formally  inverting the
class $S$ of cohomology isomorphisms, to form
$\D_{\fc}^{\fc'}(\ca)=S^{-1}\cc$,
The information retained isn't much more than the cohomology of
the complex. Ordinary homological algebra teaches us that there are
really only two things you can do with cohomology:
\be
\item
Shift the degrees.
\item
Form the the long exact sequence in cohomology that comes from a short
exact sequence of cochain complexes.
\ee
The structure of a \emph{triangulated category,} formalized in
Definition~\ref{D1.5}~(i) and (ii) below, encapsulates this:
Definition~\ref{D1.5}(i) gives the
shifting of degrees, while Definition~\ref{D1.5}(ii) is the
abstract version
of the long exact sequence in cohomology coming from a short exact
sequence of cochain complexes. See Example~\ref{E1.7} for more detail:
we spell out the recipe that endows $\D_{\fc}^{\fc'}(\ca)$
with the structure of a triangulated
category, and do so by steering as close as possible to the simple,
motivating idea.
\ermk

\dfn{D1.5}
To give the additive
category $\ct$ the structure of a \emph{triangulated category}
we must:
\be
\item
Specify an invertible additive endofunctor $\ct\la\ct$. In this article
we will denote it $[1]$ and have it act on the right:
thus it takes the object $X$ and the morphism
$f$ in $\ct$ to $X[1]$ and $f[1]$, respectively.
\setcounter{enumiv}{\value{enumi}}
\ee
Before we continue the definition we set up
\edfn

\ntn{N1.5973}
For the purpose of the current definition (Definition~\ref{D1.5}) we 
adopt the following convention. With $[1]:\ct\la\ct$ the endofunctor
of Definition~\ref{D1.5}(i), a 
\emph{candidate triangle} is any
three composable morphisms
$X\stackrel f\la Y\stackrel g\la Z\stackrel h\la X[1]$ in the
category $\ct$. The candidate triangles form a category, a
\emph{morphism} of candidate triangles 
is defined to be a commutative diagram in $\ct$
\[\xymatrix@C+20pt{
X\ar[r]^-f\ar[d]^-u & Y\ar[r]^-g\ar[d]^-v & Z\ar[r]^-h\ar[d]^-w & X[1]\ar[d]^-{u[1]} \\
X'\ar[r]^-{f'} & Y'\ar[r]^-{g'} & Z'\ar[r]^-{h'} & X'[1] 
}\]
which we view as a morphism from the top to the bottom row.
The composition of morphisms of candidate triangles is the obvious.
\entn

\nin
{\bf Continuation of Definition~\ref{D1.5}, now that the notation has been
explained.}\ \ 
In addition to the invertible endomorphism $[1]:\ct\la\ct$ of (i)
we must specify
\emph{
\be
\setcounter{enumi}{\value{enumiv}}
\item
A full subcategory of the category of candidate triangles,
whose objects will be called \emph{exact triangles.} [In some parts 
of the literature they go by the name \emph{distinguished triangles.}]
\setcounter{enumiv}{\value{enumi}}
\ee
For $\ct$ to qualify as a triangulated category the data
of (i) and (ii) above must satisfy the following axioms:
\begin{description}
\item[{[TR1]}]
Any candidate triangle isomorphic 
to an exact triangle is an exact triangle. For any object
$X\in\ct$ the diagram $0\la X\stackrel\id\la X\la 0$ is an exact triangle.  
Any morphism $f:X\la Y$ may be completed to an exact triangle
$X\stackrel f\la Y\stackrel g\la Z\stackrel h\la X[1]$.
\item[{[TR2]}]
Any rotation of an exact triangle is exact. That is: 
$X\stackrel f\la Y\stackrel g\la Z\stackrel h\la X[1]$
is an exact triangle if and only if
$Y\stackrel {-g}\la Z\stackrel {-h}\la X[1]\stackrel {-f[1]}\la Y[1]$
is.
\item[{[TR3+4]}]
Given a commutative diagram, where the rows are exact triangles,
\[\xymatrix@C+20pt{
X\ar[r]^-f\ar[d]^-u & Y\ar[r]^-g\ar[d]^-v & Z\ar[r]^-h & X[1] \\
X'\ar[r]^-{f'} & Y'\ar[r]^-{g'} & Z'\ar[r]^-{h'} & X'[1] 
}\]
we may complete it to a morphism of exact triangles, that is
a commutative diagram
\[\xymatrix@C+20pt{
X\ar[r]^-f\ar[d]^-u & Y\ar[r]^-g\ar[d]^-v & Z\ar[r]^-h\ar[d]^-w & X[1]\ar[d]^-{u[1]} \\
X'\ar[r]^-{f'} & Y'\ar[r]^-{g'} & Z'\ar[r]^-{h'} & X'[1] 
}\]
Moreover: we can do it in such a way that
\[\xymatrix@C+40pt{
  Y\oplus X'\ar[r]^-{\left(\begin{array}{rr}
      -g & 0\\
      v & f'
    \end{array}\right)} &
  Z\oplus Y'\ar[r]^-{\left(\begin{array}{rr}
      -h & 0\\
      w & g'
    \end{array}\right)} &
  X[1]\oplus Z'\ar[r]^-{\left(\begin{array}{rr}
      -f[1] & 0\\
      u[1] & h'
    \end{array}\right)}   &
  Y[1]\oplus X'[1] 
}\]
is an exact triangle.
\end{description}
%\end{emph}
}
\medskip

\exm{E1.7}
We have asserted that the category $\D_\fc^{\fc'}(\ca)$
is triangulated. It is only
fair to tell the reader what is the endofunctor
$[1]:\D_\fc^{\fc'}(\ca)\la\D_\fc^{\fc'}(\ca)$
and what are the
exact triangles. The endofunctor $[1]$, called the shift or suspension,
is easy: it takes
the cochain complex $A^*$, that is the diagram
\[\xymatrix@C+20pt{
\cdots \ar[r]& A^{-2}  \ar[r]^-{\partial^{-2}}& A^{-1}  \ar[r]^-{\partial^{-1}}& A^{0}  \ar[r]^-{\partial^{0}}
&  A^{1} \ar[r]^-{\partial^{1}}& A^{2} \ar[r] &\cdots
}\]
to the cochain complex $\big(A[1]\big)^*$ below:
\[\xymatrix@C+20pt{
\cdots \ar[r]& A^{-1}  \ar[r]^-{-\partial^{-1}}& A^{0}  \ar[r]^-{-\partial^{0}}
&  A^{1} \ar[r]^-{-\partial^{1}}& A^{2} \ar[r]^-{-\partial^{2}} & A^{3}  \ar[r]&\cdots
}\]
In words: we shift the complex to the left by one, that is
$\big(A[1]\big)^n=A^{n+1}$, and the maps all change signs.
This deals with objects.

If $f^*:A^*\la B^*$ is a cochain map
\[\xymatrix@C+10pt{
\cdots \ar[r]& A^{-2}  \ar[d]^-{f^{-2}}\ar[r]^-{\partial_A^{-2}}& A^{-1} \ar[d]^-{f^{-1}} \ar[r]^-{\partial_A^{-1}}& A^{0}\ar[d]^-{f^{0}}  \ar[r]^-{\partial_A^{0}}
&  A^{1} \ar[d]^-{f^{1}}\ar[r]^-{\partial_A^{1}}& A^{2}\ar[d]^-{f^{2}} \ar[r] &\cdots\\
\cdots \ar[r]& B^{-2}  \ar[r]_-{\partial_B^{-2}}& B^{-1}  \ar[r]_-{\partial_B^{-1}}& B^{0}  \ar[r]_-{\partial_B^{0}}
&  B^{1} \ar[r]_-{\partial_B^{1}}& B^{2} \ar[r] &\cdots
}\]
then $\big(f[1]\big)^*$ is the cochain map
\[\xymatrix@C+10pt{
\cdots \ar[r]& A^{-1} \ar[d]^-{f^{-1}} \ar[r]^-{-\partial_A^{-1}}& A^{0}\ar[d]^-{f^{0}}  \ar[r]^-{-\partial_A^{0}}
&  A^{1} \ar[d]^-{f^{1}}\ar[r]^-{-\partial_A^{1}}& A^{2}\ar[d]^-{f^{2}}\ar[r]^-{-\partial_A^{2}}& \ar[r]
 A^{3}  \ar[d]^-{f^{3}}&\cdots\\
\cdots \ar[r]&  B^{-1}  \ar[r]_-{-\partial_B^{-1}}& B^{0}  \ar[r]_-{-\partial_B^{0}}
&  B^{1} \ar[r]_-{-\partial_B^{1}}& B^{2}\ar[r]_-{-\partial_B^{2}}& B^{3}  \ar[r] &\cdots
}\]
This defines what the functor $[1]$ does to cochain maps, and we extend
to arbitrary morphisms in $\D_\fc^{\fc'}(\ca)$ by the
universal property of the localization
process.

To spell this out a bit, as in Explanation~\ref{X1.3}: let $\cc$ be the
category with the same objects as $\D_\fc^{\fc'}(\ca)$
but where the morphisms are
the cochain maps. We have defined a functor $[1]:\cc\la\cc$, which takes
the class $S\subset\text{Mor}(\cc)$ of cohomology isomorphisms
to itself. The composite $\cc\stackrel{[1]}\la\cc\stackrel F\la S^{-1}\cc$
is a functor from $\cc$ to the category $S^{-1}\cc=\D_\fc^{\fc'}(\ca)$,
which takes
the morphisms in $S$ to isomorphisms. By the universal property it
factors uniquely through $F$, that is there exists a commutative square
\[\xymatrix@C+10pt{
\cc \ar[r]^-{[1]} \ar[d]_F & \cc\ar[d]^F \\
S^{-1}\cc \ar@{.>}[r]^{\exists!} &S^{-1}\cc
}\]
We declare $[1]:\D_\fc^{\fc'}(\ca)\la\D_\fc^{\fc'}(\ca)$ to be the
unique map $S^{-1}\cc\la S^{-1}\cc$
making the square commute.

It remains to describe the exact triangles---Remark~\ref{R1.4} provided the
intuition, it told us that
the exact triangles should be the formalization of the long exact
sequence in cohomology coming from a short exact sequence
of cochain complexes.
We propose to give the skeleton of the construction below, and
the reader interested in more detail is referred
to the appendices.

Suppose therefore that we are given a commutative
diagram in $\ca$
\[\xymatrix@C+10pt{
\cdots \ar[r]& X^{-2}  \ar[d]\ar[r]& X^{-1} \ar[d] \ar[r]& X^{0}\ar[d]  \ar[r]
&  X^{1} \ar[d]\ar[r]& X^{2}\ar[d] \ar[r] &\cdots\\
\cdots \ar[r]& Y^{-2}  \ar[d]\ar[r]& Y^{-1} \ar[d] \ar[r]& Y^{0}\ar[d]  \ar[r]
&  Y^{1} \ar[d]\ar[r]& Y^{2}\ar[d] \ar[r] &\cdots\\
\cdots \ar[r]& Z^{-2}  \ar[r]& Z^{-1}  \ar[r]& Z^{0}  \ar[r]
&  Z^{1} \ar[r]& Z^{2} \ar[r] &\cdots
}\]
where the rows are objects of $\D_\fc^{\fc'}(\ca)$,
that is cochain complexes satisfying the hypotheses. So far
we may view the above as morphisms
$X^*\stackrel{f^*}\la Y^*\stackrel{g^*}\la Z^*$
in the category $\D_\fc^{\fc'}(\ca)$. Assume
further that, for each $i\in\zz$, the sequence
$X^i\stackrel{f^i}\la Y^i\stackrel{g^i}\la Z^i$ is
split exact---it's easier to deal with degreewise
split short exact sequences, in Appendix~\ref{A2} the reader will
see that,
up to isomorphism 
in $\D(R)$, this suffices. We next  want to mimic the process that
produces the differential of the long exact sequence in
cohomology. Choose, for each $i\in\zz$, a splitting $\theta^i:Z^i\la Y^i$
of the map $g^i:Y^i\la Z^i$. Now for each $i$ we have the diagram
\[\xymatrix@C+20pt{
Z^i\ar[r]^-{\theta^i}\ar[d]_-{\partial_Z^i} & Y^i\ar[d]^-{\partial_Y^i}
\ar[r]^-{g^i}  & Z^i\ar[d]^-{\partial_Z^i}\\
Z^{i+1}\ar[r]^-{\theta^{i+1}} & Y^{i+1}\ar[r]^-{g^{i+1}} &Z^{i+1}
}\]
If we delete the middle column the resulting square commutes---the composites
of the horizontal maps are identities. If we delete the left column
the square is commutative because it is part of the diagram defining
the cochain map $g^*$. It follows that, in the diagram below,
\[\xymatrix@C+20pt{
Z^i\ar[r]^-{\theta^i}\ar[d]_-{\partial_Z^i} & Y^i\ar[d]^-{\partial_Y^i} \\
Z^{i+1}\ar[r]^-{\theta^{i+1}} & Y^{i+1}\ar[r]^-{g^{i+1}} &Z^{i+1}
}\]
the
two composites from top left to bottom right are equal.
Thus the difference
$\theta^{i+1}\partial_Z^i-\partial_Y^i\theta^i$ is
annihilated by the map $g^{i+1}:Y^{i+1}\la Z^{i+1}$, hence
$\theta^{i+1}\partial_Z^i-\partial_Y^i\theta^i$ must factor uniquely
through the kernel of $g^{i+1}$, it can be written uniquely 
as the composite
$Z^i\stackrel{h^{i}}\la X^{i+1}\stackrel{f^{i+1}}\la Y^{i+1}$.
In Appendix~\ref{A1}, the reader can see that the following is a cochain map
\[\xymatrix@C+10pt{
\cdots \ar[r]& Z^{-2}  \ar[d]^-{h^{-2}}\ar[r]^-{\partial_Z^{-2}}& Z^{-1} \ar[d]^-{h^{-1}} \ar[r]^-{\partial_Z^{-1}}& Z^{0}\ar[d]^-{h^{0}}  \ar[r]^-{\partial_Z^{-0}}
&  Z^{1} \ar[d]^-{h^{1}}\ar[r]^-{\partial_Z^{1}}& Z^{2}\ar[d]^--{h^{2}} \ar[r] &\cdots\\
\cdots \ar[r]& X^{-1}  \ar[r]_-{-\partial_X^{-1}}& X^{0}  \ar[r]_-{-\partial_X^{0}}& X^{1}  \ar[r]_-{-\partial_X^{1}}
&  X^{2} \ar[r]_-{-\partial_X^{2}}& X^{3} \ar[r] &\cdots
}\]
Thus we have constructed in the category $\D_\fc^{\fc'}(\ca)$ a diagram
$X^*\stackrel {f^*}\la Y^*\stackrel {g^*}\la Z^*\stackrel {h^*}\la X^*[1]$.
We declare
\be
\item
The exact triangles in $\D_\fc^{\fc'}(\ca)$ are all the isomorphs,
in $\D_\fc^{\fc'}(\ca)$, of diagrams
that come from our construction.
\ee
It needs to be checked that [TR1], [TR2] and [TR3+4] of Definition~\ref{D1.5}
are satisfied, the
reader can amuse herself with this.
\eexm

For future reference we recall:

\ntn{N1.8}
If $\ct$ is a triangulated category and $n\in\zz$ is an integer, then
$[n]$ will be our shorthand for the endofunctor $[1]^n:\ct\la\ct$.
Also: we will often lazily abbreviate ``exact triangle'' to just ``triangle''.
\entn

\dfn{D1.9}
A full subcategory $\cs\subset\ct$ is called \emph{triangulated} if
$0\in\cs$, if $\cs[1]=\cs$, and if, whenever $X,Y\in\cs$ and
there exists in $\ct$ a triangle
$X\la Y\la Z\la X[1]$, we must also have $Z\in\cs$. The subcategory $\cs$
is \emph{thick} if it is triangulated, as well as closed in $\ct$
under direct summands.
\edfn

Now that we have recalled the notion of triangulated categories,
as well as thick and triangulated subcategories, it is time
to remember the other two building blocks of the theory we plan to
introduce: compact generators and {\it t}--structures. We begin with

\dfn{D1.11}
Let $\ct$ be a triangulated category with coproducts. An object $C\in\ct$
is \emph{compact} if the functor $\Hom(C,-)$ respects coproducts. A set
of compact
objects $\{G_i,\,i\in I\}$ is said to \emph{generate} the category $\ct$
if the following equivalent conditions hold
\be
\item
If $X\in\ct$ is an object, and if
$\Hom\big(G_i,X[n]\big)\cong0$ for all $i\in I$ and
all $n\in\zz$, then $X\cong0$.
\item
If a triangulated subcategory $\cs\subset\ct$ is closed under coproducts
and contains the objects $\{G_i,\,i\in I\}$, then $\cs=\ct$.
\ee
If the category $\ct$ contains a set of compact generators it is
called \emph{compactly generated.}
\edfn

\rmk{R1.13}
The equivalence of (i) and (ii) in Definition~\ref{D1.11} is not meant to be
obvious, but it is a standard result. We will mostly be interested in the
situation where the category $\ct$ is compactly generated and, moreover,
the set of compact generators may be chosen to consist of
a single element. That is: for some compact object $G\in\ct$ the set
$\{G\}$ generates, as in Definition~\ref{D1.11}~(i) or (ii).
\ermk

\exm{E1.15}
The categories $\D(R)$  and $\Dqc(X)$ of
Example~\ref{E1.100372} (i) and (ii) both have coproducts\footnote{
  In this example it helps to know the calulus of fractions of
  $\D_\fc^{\fc'}(\ca)=S^{-1}\cc$
 mentioned in Example~\ref{X1.3}.
 After all we are making assertions about morphisms in $\D_\fc^{\fc'}(\ca)$: to say
 that an object is a coproduct is a universal property for certain morphisms.
Moreover we also make an assertion about $\Hom_{\D(R)}^{}(R,-)$.}:
the coproduct of a family of cochain
complexes
\[\xymatrix@C+20pt{
\cdots \ar[r]& A_\lambda^{-2}  \ar[r]& A_\lambda^{-1}  \ar[r]& A_\lambda^{0}  \ar[r]
&  A_\lambda^{1} \ar[r]& A_\lambda^{2} \ar[r] &\cdots
}\]
turns out to be nothing other than
\[\xymatrix@C+2pt{
\cdots \ar[r]& \ds\coprod_{\lambda\in\Lambda}A_\lambda^{-2}  \ar[r]& \ds\coprod_{\lambda\in\Lambda}A_\lambda^{-1}  \ar[r]& \ds\coprod_{\lambda\in\Lambda}A_\lambda^{0}  \ar[r]
&  \ds\coprod_{\lambda\in\Lambda}A_\lambda^{1} \ar[r]& \ds\coprod_{\lambda\in\Lambda}A_\lambda^{2} \ar[r] &\cdots
}\]
It's clear that the formula above does not work for
the categories $\dperf X$ and $\dcoh(X)$ of
Example~\ref{E1.100372} (iii) and (iv), if we take a giant direct sum
of complexes satisfying the restrictions the resulting complex
will fail to satisfy the restrictions. And it's not just that the
formulas don't work, the categories $\dperf X$ and $\dcoh(X)$ don't
have coproducts.

Now for compact generators.
If $R\in\D(R)$ stands for the cochain complex
\[\xymatrix@C+20pt{
\cdots \ar[r]& 0  \ar[r]& 0  \ar[r]&R  \ar[r]
&  0 \ar[r]& 0\ar[r] &\cdots
}\]
that is the complex whose only nonzero entry is the module $R$ in degree 0,
then it can be shown that there is
an isomorphism of functors $\Hom_{\D(R)}^{}(R,-)\cong H^0(-)$.
The functor $H^0(-)$ obviously
respects coproducts, hence so does $\Hom_{\D(R)}^{}(R,-)$;
that is the object $R\in\D(R)$ is compact.

Next observe that, if $X\in\D(R)$ is an object such that
$H^n(X)\cong\Hom\big(R,X[n]\big)\cong0$ for all $n\in\zz$, then $X$ is acyclic;
its cohomology all vanishes. The cochain map $0\la X$ is an isomorphism
in cohomology, hence an isomorphism in $\D(R)$. That is: $X\cong0$. Thus
the compact object $R\in\D(R)$ satsifies Definition~\ref{D1.11}(i), it
is a compact generator.

The category $\D(R)$ is compactly generated, and more precisely we have learned
that the object $R\in\D(R)$ is a single compact generator.

Not so easy is the fact that, if $X$ is a quasicompact,
quasiseparated scheme, then the category $\Dqc(X)$ also has a
single compact generator. This is a theorem,
proved in Bondal and Van den
Bergh~\cite[Theorem~3.1.1(ii)]{BondalvandenBergh04}.
\eexm

\ntn{N1.107589}
Let $\ct$ be a triangulated category with coproducts. It is standard
to denote by $\ct^c$ the full subcategory, whose objects are the
compact objects in $\ct$. It isn't difficult to show that $\ct^c$ is
always a thick subcategory of $\ct$, as in Definition~\ref{D1.9}.
In the case where
$\ct=\Dqc(X)$ the category $\ct^c$ turns out to be the 
$\dperf X$ of Example~\ref{E1.100372}(iii),
the reader can find this fact in
Bondal and Van den
Bergh~\cite[Theorem~3.1.1(i)]{BondalvandenBergh04}.
\entn

We also need to recall {\it t}--structures, and we plan to
begin with the concrete. But first a reminder.

\rmd{R1.90909}
Let $\ct$ be an triangulated category, and $\ca$ an abelian category.
A functor $H:\ct\la \ca$ is called \emph{homological} if, for
every triangle $A\la B\la C\la A[1]$ in $\ct$, the sequence
$H(A)\la H(B)\la H(C)$ is exact in $\ca$.

From the axiom [TR2] of Definition~\ref{D1.5}---the axiom telling us
that any rotation of a triangle is a triangle---it follows that the functor
$H$ must take 
a triangle in $\ct$ to a long exact
sequence.
\ermd

\exm{E1.59503}
It follows from the axioms of triangulated categories that
all representable functors are homological. That is: if $\ct$ is
a triangulated category and $A\in\ct$ is an object, then
$\Hom(A,-)$ and $\Hom(-,A)$ are, respectively, homological
functors $\ct\la\ab$ and $\ct\op\la\ab$, where $\ab$ is the category
of abelian groups. 

The functor $H:\D(R)\la\Mod R$, taking a complex to its zeroth cohomology,
is homological. In Example~\ref{E1.15} we were told that $H(-)\cong\Hom(R,-)$,
that is the functor $H$ is a special case of the previous paragraph,
it is a representable functor.

On the categories $\Dqc(X)$, $\dperf X$ and $\dcoh(X)$ the homological
functor we will usually consider is traditionally denoted
$\ch$, and takes its values in the abelian category $\Mod{\co_X^{}}$
of sheaves of $\co_X^{}$--modules. Again: the functor $\ch$ just
takes a complex of sheaves to the zeroth cohomology
sheaf.\footnote{For
  the non-algebraic-geometers: the letter $H$ is taken, it usually means
  another homological functor.}

The fact that $H$ and $\ch$ are homological is by the construction
of triangles,
see Example~\ref{E1.7}---it comes down to the statement
that the long exact sequence coming from a short exact sequence of
cochain complexes is exact.
\eexm

And finally we turn to {\it t}--structures, introducing them by
example.

\exm{E1.17}
In the category $\ct=\D(R)$
we define two full subcategories by the formula
\be 
\item $\ct^{\leq0}\quad=\quad\{A\in\D(R)\mid H\big(A[i]\big)=0\text{ for all }i>0\}$
\item $\ct^{\geq0}\quad=\quad\{A\in\D(R)\mid H\big(A[i]\big)=0\text{ for all }i<0\}$
\setcounter{enumiv}{\value{enumi}}
\ee
While in the case where
$\ct$ is either of the categories $\Dqc(X)$ or $\dcoh(X)$, the formula is
\be 
\setcounter{enumi}{\value{enumiv}}
\item $\ct^{\leq0}\quad=\quad\{A\in\ct\mid \ch\big(A[i]\big)=0\text{ for all }i>0\}$
\item $\ct^{\geq0}\quad=\quad\{A\in\ct\mid \ch\big(A[i]\big)=0\text{ for all }i<0\}$
\setcounter{enumiv}{\value{enumi}}
\ee
These pairs of subcategories, in each of $\D(R)$,
$\Dqc(X)$ and $\dcoh(X)$, define a \tstr.
For each of the three categories the particular \tstr\ above is
traditionally called
the
\emph{standard \tstr.} The category $\dperf X$ does not usually have
a nontrivial \tstr.

Let us next give the formal definition:
\eexm

\dfn{D1.19}
A \emph{$t$--structure} on a triangulated category $\ct$ is a pair
of full subcategories $\big(\ct^{\leq0}_{},\ct^{\geq0}_{}\big)$ satisfying 
\be 
\item   \(
\ct^{\leq0}_{}[1]\subset \ct^{\leq0}_{}\qquad\text{ and }\qquad
\ct^{\geq0}_{}\subset \ct^{\geq0}_{}[1]
\)
\item   \(
\Hom\Big(\ct^{\leq0}_{}[1]\,\,,\,\,
\ct^{\geq0}_{}\Big)=0
\)
\item
Every object $B\in\ct$ admits a triangle $A\la B\la C\la$ with
$A\in\ct^{\leq0}_{}[1]$ and $C\in\ct^{\geq0}_{}$.
\ee 
\edfn

\rmk{R1.21}
It can be checked that the pairs of subcategories
of Example~\ref{E1.17} satisfy parts (i), (ii) and (iii) of
Definition~\ref{D1.19}, they do provide {\it t}--structures on
each of
$\D(R)$, $\Dqc(X)$ and $\dcoh(X)$.
\ermk

We have now introduced all the players: triangulated categories, compact
generators and {\it t}--structures. We end the section recalling
certain
standard shorthand 
conventions.

\ntn{N1.23}
Let $\ct$ be a triangulated category with a \tstr\ 
$\big(\ct^{\leq0}_{},\ct^{\geq0}_{}\big)$. Then
\be
\item
For any integer $n\in\zz$ we set
\[
\ct^{\leq n}_{}=\ct^{\leq 0}_{}[-n]\qquad{ and }\qquad
\ct^{\geq n}_{}=\ct^{\geq 0}_{}[-n]
\]
\item Furthermore, we adopt the conventions
  \[
  \ct^-=\bigcup_{n\in\nn}\ct^{\leq n},\qquad
  \ct^+=\bigcup_{n\in\nn}\ct^{\geq-n},\qquad
  \ct^b=\ct^-\cap\ct^+.
  \]
\ee
\entn

\section{Approximability---the intuition, which comes from $\D(R)$}
\label{S2}

In the last section we recalled, for the benefit of the non-expert, some
standard facts about triangulated categories, compact generators and
{\it t}--structures---as well as the special cases
that play a big role in this article, namely $\D(R)$, $\Dqc(X)$, $\dperf X$
and $\dcoh(X)$. It's time to move on to the subject matter
of this article: approximability. As we've tried to do throughout, we
will proceed from the concrete to the abstract. Let us
therefore first study what it all means for the category $\D(R)$, when $R$ is a
ring.

The category $\D(R)$ has a standard \tstr,
see Example~\ref{E1.17}, Definition~\ref{D1.19} and Remark~\ref{R1.21}.
Suppose we are given an object $F^*\in\D(R)^{\leq0}$, meaning
a cochain complex
\[\xymatrix@C+10pt{
\cdots \ar[r]& F^{-2}  \ar[r]& F^{-1}  \ar[r]& F^{0}  \ar[r]
&  F^{1} \ar[r]& F^{2} \ar[r] &\cdots
}\]
such that $H^i(F^*)=0$ for all $i>0$.
Then $F^*$ has a projective resolution. We can produce a
cochain map
\[\xymatrix@C+10pt{
\cdots \ar[r]& P^{-2}  \ar[d]\ar[r]& P^{-1} \ar[d] \ar[r]& P^{0}\ar[d]  \ar[r]
&  0 \ar[d]\ar[r]& 0\ar[d] \ar[r] &\cdots\\
\cdots \ar[r]& F^{-2}  \ar[r]& F^{-1}  \ar[r]& F^{0}  \ar[r]
&  F^{1} \ar[r]& F^{2} \ar[r] &\cdots
}\]
inducing an isomorphism in cohomology, and so that each $P^i$ is a
projective $R$--module.
This gives us, in the category $\D(R)$, an isomorphism $P^*\la F^*$.
Now consider
\[\xymatrix@C+2pt{
\cdots \ar[r]& 0  \ar[d]\ar[r]& P^{-n}  \ar[d]\ar[r]& \cdots  \ar[r]& P^{-1}\ar[d]  \ar[r]
&   P^{0} \ar[d]\ar[r]& 0\ar[d] \ar[r] &\cdots\\
\cdots \ar[r]& P^{-n-1}  \ar[d]\ar[r]& P^{-n}  \ar[d]\ar[r]& \cdots  \ar[r]& P^{-1}\ar[d]  \ar[r]
&   P^{0} \ar[d]\ar[r]& 0\ar[d] \ar[r] &\cdots\\
\cdots \ar[r]& P^{-n-1}  \ar[r]& 0 \ar[r]& \cdots  \ar[r]& 0  \ar[r]
&   0 \ar[r]& 0 \ar[r] &\cdots
}\]
This yields a pair of cochain maps
$E^*_n\stackrel{f_n^*}\la P^*\stackrel{g_n^*}\la D^*_n$ so that,
in each degree $i$, the maps
$E^i_n\stackrel{f_n^i}\la P^i\stackrel{g_n^i}\la D^i_n$
deliver a split exact sequence of $R$--modules. Example~\ref{E1.7}
constructs for us a cochain map
$h^*_n:D_n^*\la E^*_n[1]$ so that the diagram
$E^*_n\stackrel{f_n^*}\la P^*\stackrel{g_n^*}\la D^*_n\stackrel{h_n^*}\la E^*_n[1]$
is an exact triangle. The isomorphism $P^*\la F^*$ in the category
$\D(R)$, coupled with the fact that any isomorph of a
triangle is a
triangle, produces in $\D(R)$ a triangle which we will write
$E^*_n\stackrel{f_n}\la F^*\stackrel{g_n}\la D^*_n\stackrel{h^*_n}\la E^*_n[1]$.

\smr{S2.1}
Given an object $F^*\in\D(R)^{\leq0}$ and an integer $n\geq0$
we have constructed,
in $\D(R)$, a triangle
$E^*_n\stackrel{f_n}\la F^*\stackrel{g_n}\la D^*_n\stackrel{h_n^*}\la E^*_n[1]$.
This triangle is such that
$D^*_n\in\D(R)^{\leq-n-1}$, while $E^*_n$ is not too complicated.

In the Introduction we mentioned that we will view the objects
$D^*_n$ as ``small'' with respect to the metric induced by the \tstr.
Up to an arbitrarily small ``correction term'' $D^*_n$, we
have a way of approximating
the object $F^*$ by the object $E^*_n$ which we
view as simpler. In order to formalize the idea
we need to make precise what we mean by
saying that $E^*_n$ is ``not too complicated''.
We will do this in the next section.
\esmr

\section{Measuring the complexity of an object}
\label{S3}

As we have said in the Introduction, measuring how complicated 
an object is will involve a small tweak of an idea from 
Bondal and Van~den~Bergh~\cite{BondalvandenBergh04}.
We remind the reader.

\rmd{R3.1}
Let $\ct$ be a triangulated category, possibly with coproducts,
and let $\ca, \cb\subset\ct$
be full subcategories.  We define the full subcategories
\be 
\item 
$
\ds\ca*\cb\quad=\quad\left\{x\in\ct\left|\begin{array}{c}
\text{ there exists in $\ct$ a triangle }a\la x\la b\\
\text{ with }a\in\ca\text{ and }\,b\in\cb\end{array}\right.\right\}
$
\item
$\text{add}(\ca)$: this consists of all finite coproducts of objects of $\ca$.
\item  
Assume $\ct$ has coproducts. Define
$\text{Add}(\ca)$ to consist of all coproducts of objects of $\ca$.
\item
  $\text{smd}(\ca)$: the category of
  all direct summands of objects of $\ca$.
\ee
\ermd  

\rmk{R3.3}
Reminder~\ref{R3.1}(i) is as in \cite[1.3.9]{BeiBerDel82}, while 
Reminder~\ref{R3.1}(iv) is identical with \cite[beginning of
2.2]{BondalvandenBergh04}.
Reminder~\ref{R3.1}~(ii) and (iii) follow the usual conventions in
representation
theory; in \cite[beginning of
2.2]{BondalvandenBergh04} the authors adopt the
(unconventional) notation that $\add(\ca)$ and $\Add(\ca)$
are also closed under the
suspension---thus $\add(\ca)$ as defined in
\cite{BondalvandenBergh04} is what we would denote
$\add\big(\bigcup_{n=-\infty}^\infty\ca[n]\big)$. 
The definitions that follow are therefore slightly
different from \cite{BondalvandenBergh04}, and it is this small
tweak that makes all the difference---with the tweaked definitions,
approximations turn out to exist in great generality.

And now we come to the key definition: we're about to measure
how much effort goes into constructing an object $X$ out of
some given full subcategory $\ca\subset\ct$. In practice
our usual choice for $\ca$ will be a $\ca=\{G\}$, the subcategory
with just a single object $G$, which we will often assume
to be a compact generator.
\ermk

\dfn{D3.5}
Let $\ct$ be a triangulated category, possibly with coproducts,
let $\ca\subset\ct$
be a full subcategory and let 
$m\leq n$ be integers, possibly infinite.
We define the full subcategories
\be 
\item
$\ca[m,n]=\cup_{i=m}^n\ca[-i]$.
\item 
$\genu\ca1{m,n}=\Smr\Big[\add\big(\ca[m,n]\big)\Big]$.
\item
$\ogenu\ca1{m,n}=\Smr\Big[\Add\big(\ca[m,n]\big)\Big]$.\ \ \ 
 \emph{[This definition assumes $\ct$ has coproducts].}
\setcounter{enumiv}{\value{enumi}}
\ee
Now let $\ell>0$ be an integer,
and assume the categories  $\genu\ca k{m,n}$ and
$\ogenu\ca k{m,n}$ have been defined for $k$ in the range $1\leq k\leq\ell$.
We proceed inductively to set 
\be
\setcounter{enumi}{\value{enumiv}}
\item
$\genu\ca{\ell+1}{m,n}=\Smr\Big[\genu\ca1{m,n}*\genu\ca\ell{m,n}\Big]$.
\item   
$\ogenu\ca{\ell+1}{m,n}=\Smr\Big[\ogenu\ca1{m,n}*\ogenu\ca\ell{m,n}\Big]$.
\ \ \ \emph{[This definition assumes $\ct$ has coproducts].}
\ee 
\edfn

\exm{E3.7}
Let us go back to our favorite example $\D(R)$. Suppose
$\ca=\{R\}$ is the category with a single object $R$, and we
will now proceed
to say something about the subcategories $\genu R\ell{-n,0}\subset\D(R)$.
Let us start with
\be 
\item  
$\genu R{1}{-n,0}$: this turns out to be the category of all 
isomorphs in $\D(R)$ of the cochain complexes 
\[\xymatrix@C-2pt{
\cdots \ar[r]& 0  \ar[r]& P^{-n}  \ar[r]^-0& \cdots  \ar[r]^-0& P^{-1}  \ar[r]^-0
&   P^{0} \ar[r]& 0\ar[r] &\cdots
}\]
with $P^i$ finitely generated and projective.
\setcounter{enumiv}{\value{enumi}}
\ee
This much is basically true by construction. We start with the object $R$, and
in Definition~\ref{D3.5}(i) we form the 
category $R[-n,0]=\{R[i],\,0\leq i\leq n\}$ with finitely many objects.
And then Definition~\ref{D3.5}(ii) allows us to first form finite coproducts
of objects in $R[-n,0]$, meaning cochain complexes
\[\xymatrix@C-2pt{
\cdots \ar[r]& 0  \ar[r]& P^{-n}  \ar[r]^-0& \cdots  \ar[r]^-0& P^{-1}  \ar[r]^-0
&   P^{0} \ar[r]& 0\ar[r] &\cdots
}\]
with each $P^i$ a finitely generated free module, and then we are permitted
direct summands in $\D(R)$ of the above. It may be shown that
these are all isomorphic to complexes as above, but where we allow the $P^i$
to be finitely generated and projective.

This was the easy part. Now the categories 
$\genu R\ell{-n,0}$ grow as $\ell$ grows, but it's a little unclear
how fast.
They all contain $\genu R{1}{-n,0}$, and are all contained in the subcategory
$\cs\subset\D(R)$ of objects isomorphic in $\D(R)$ to cochain complexes
\[\xymatrix@C-2pt{
\cdots \ar[r]& 0  \ar[r]& P^{-n}  \ar[r]& \cdots  \ar[r]& P^{-1}  \ar[r]
&   P^{0} \ar[r]& 0 \ar[r] &\cdots
}\]
with $P^i$ finitely generated and projective. Unlike the
category $\genu R1{-n,0}$ of (i) above, for
objects in $\cs$ the maps
$P^i\la P^{i+1}$
are unconstrained---beyond (of course) the standing assumption that
all composites $P^i\la P^{i+1}\la P^{i+2}$ must vanish, the objects of
$\cs\subset\D(R)$ must be cochain complexes.

What turns out to be true is 
\be
\setcounter{enumi}{\value{enumiv}}
\item
$\genu R{n+1}{-n,0}=\cs$; hence $\genu R\ell{-n,0}=\cs$ for all $\ell\geq n+1$.
\ee
We leave to the reader the proofs of the
assertions made in this Example\footnote{The proofs are easy for the reader
  familiar with the calculus of fractions of Explanation~\ref{X1.3}. Other
  readers are asked to accept the assertions on faith.}.
\eexm

\exm{E3.9}
Let us stay with our favorite example $\D(R)$, and let us continue
to put $\ca=\{R\}$, that is
$\ca$ is the full subcategory of $\D(R)$ with the single object $R$.
We now want to work out what are the categories $\ogenu R{\ell}{-n,0}$.
The discussion turns out to be much the same as in Example~\ref{E3.7},
and the summary of the results is
\be 
\item  
The category $\ogenu R{1}{-n,0}$ consists of all isomorphs in
$\D(R)$ of cochain complexes
\[\xymatrix@C-2pt{
\cdots \ar[r]& 0  \ar[r]& P^{-n}  \ar[r]^-0& \cdots  \ar[r]^-0& P^{-1}  \ar[r]^-0
&   P^{0} \ar[r]& 0\ar[r] &\cdots
}\]
with $P^i$ projective.
\item
The category $\ogenu R{n+1}{-n,0}$ consists of all isomorphs in $\D(R)$
of cochain  complexes
\[\xymatrix@C-2pt{
\cdots \ar[r]& 0  \ar[r]& P^{-n}  \ar[r]& \cdots  \ar[r]& P^{-1}  \ar[r]
&   P^{0} \ar[r]& 0 \ar[r] &\cdots
}\]
with $P^i$ projective. Moreover: if $\ell\geq n+1$ then
$\ogenu R{n+1}{-n,0}=\ogenu R{\ell}{-n,0}$.
\ee
Thus the objects in both
$\ogenu R{n+1}{-n,0}$ and $\genu R{n+1}{-n,0}$ are isomorphic in $\D(R)$ to
complexes of projectives vanishing outside the interval $[-n,0]$, and the
difference is that  in $\ogenu R{n+1}{-n,0}$ the
projective modules are not
constrained to be finitely generated.
\eexm

\cnc{C3.11}
In the new notation we have introduced, Summary~\ref{S2.1} and
Example~\ref{E3.9} combine to say: for any object $F\in\D(R)^{\leq0}$
and any integer $n\geq0$ there exists a triangle
\[\xymatrix{
E_n\ar[r]^-f & F\ar[r]^-g &D_n\ar[r]^-h & E_n[1]
}\]
with $D_n\in\D(R)^{\leq-n-1}$ and $E_n\in\ogenu R{n+1}{-n,0}$.
\ecnc

\rmk{R3.13}
Let $\D(\proj R)^{\leq0}\subset\D(R)$ be the full subcategory,
whose objects are the isomorphs in $\D(R)$ of cochain complexes
\[\xymatrix@C-2pt{
\cdots \ar[r]& P^{-n-1}  \ar[r]& P^{-n}  \ar[r]& \cdots  \ar[r]& P^{-1}  \ar[r]
&   P^{0} \ar[r]& 0 \ar[r] &\cdots
}\]
with $P^i$ finitely generated and projective. 
Summary~\ref{S2.1} and
Example~\ref{E3.7} combine to say: for any object $F\in\D(\proj R)^{\leq0}$
and any integer $n\geq0$ there exists a triangle
\[\xymatrix{
E_n\ar[r]^-f & F\ar[r]^-g &D_n\ar[r]^-h & E_n[1]
}\]
with $D_n\in\D(R)^{\leq-n-1}$ and $E_n\in\genu R{n+1}{-n,0}$.

We will return to the category $\D(\proj R)^{\leq0}$ and to its relative
\[
\D^-(\proj R)\eq\bigcup_{n\in\nn}^{}\D(\proj R)^{\leq0}[-n]
\]
much later in the article.
\ermk

\section{The formal definition of approximability}
\label{S4}

Now that we are thoroughly prepared, approximability
becomes easy to formulate precisely:

\dfn{D4.1}
Let $\ct$ be a triangulated category with coproducts. It is
\emph{approximable} if there exists a compact generator $G\in\ct$,
a {\it t}--structure $(\ct^{\leq0},\ct^{\geq0})$, and
an integer $A>0$ so that
\be 
\item
$G[A]\in\ct^{\leq0}$ and $\Hom\big(G[-A]\,,\,\ct^{\leq0}\big)=0$.
\item
For every object $F\in\ct^{\leq0}$ there exists a triangle $E\la F\la D\la E[1]$,
with $D\in\ct^{\leq-1}$ and $E\in\ogenu G{A}{-A,A}$.
\ee 
\edfn

\exm{E4.3}
Let $R$ be a ring. In Example~\ref{E1.15} we learned that
the object $R\in\D(R)$ is a compact generator, we will take this to be
our $G$ of Definition~\ref{D4.1}. For 
the \tstr\ we choose the standard one, see Example~\ref{E1.17}.
And for our integer we set $A=1$.

It's clear that $R[1]\in \D(R)^{\leq0}$ and that
$\Hom\big(R[-1],\D(R)^{\leq0}\big)=0$. This establishes
Definition~\ref{D4.1}(i). 
Finally suppose we are given an object $F\in\D(R)^{\leq0}$.
By Conclusion~\ref{C3.11}, with $n=0$, there must exist a triangle
$E\la F\la D\la E[1]$ 
with $D\in\D(R)^{\leq-1}$ and $E\in\ogenu R1{0,0}\subset\ogenu R1{-1,1}$.
This proves that Definition~\ref{D4.1}(ii) holds.

Thus the category $\D(R)$ is approximable.
\eexm

\rmk{R4.5}
The reader might be disappointed: until now we have been stressing
that approximability will allow us to obtain arbitrarily good estimates
of the objects in any approximable category $\ct$,
and in Example~\ref{E4.3} we see that the
definition only involves a zero-order approximation.

Don't let this disturb you, it's easy to iterate and estimate the given
object $F$ to arbitrarily high order. This will manifest itself
in our theorems.
\ermk

\rmk{R4.7}
In Example~\ref{E1.100372} we told the reader that, in this survey, the
key examples of triangulated categories will be $\D(R)$, $\Dqc(X)$, $\dperf X$
and $\dcoh(X)$. The definition of approximability is tailored so that the
category $\D(R)$ is obviously approximable, as we have seen in
Example~\ref{E4.3}. What about the other three?

The categories $\dperf X$ and $\dcoh(X)$ cannot possibly be approximable,
in Example~\ref{E1.15} we learned that they don't even have coproducts.

It is a non-obvious theorem that, as long as the scheme $X$ is quasicompact
and separated, the category $\Dqc(X)$ is approximable. And it should
come as a surprise---after all approximability was modeled on the idea of
taking a projective resolution of a bounded-above
cochain complex  and then truncating, and
this is a construction that can only work in the presence
of enough projectives. There aren't enough projectives in either
the category of sheaves of $\co_X^{}$--modules, or in
its subcategory of quasicoherent sheaves.

In Example~\ref{E1.15}  we mentioned that even the  existence of a single
compact generator in $\Dqc(X)$ isn't obvious, it's a theorem of
Bondal and Van den Bergh. The existence proof isn't particularly
constructive---it doesn't give us much of a handle on this compact
generator. And the definition of approximability is the assertion that
the compact generator may be chosen to satisfy several useful
properties; it decidedly isn't clear how to prove any of them.

Given that it is going to entail
real effort to prove that $\Dqc(X)$ is approximable,
it shouldn't come as a surprise that there are far-reaching
consequences.
\ermk

And now it's time for

\section{The main theorems}
\label{S5}

\rmk{R5.1}
The theorems break up into three groups, namely
\be
\item
Theorems that produce more examples of approximable categories. So far
we have discussed in some detail the example $\D(R)$, and
then made some passing
comments about $\Dqc(X)$. See Remark~\ref{R4.7}.
\item
Formal consequences of approximability---that is structure that comes
for free, which every approximable category has.
\item
Applications to concrete examples, which teach us new
and interesting facts about old and
familiar categories.
\ee
In this section we will list the results of type (i) and (iii)
by group, doing little more than
giving formal statements. In the remainder of the article we will first
expand
on the results in group (iii), saying something about what was known before
and about the proofs, both of the new and the old versions---presumably
the reader is most likely
to be persuaded by the theory if she can see applications
that matter.

And then, towards the end of the article, we will give results in group (ii).
We hope that by then, with the reader's interest piqued by the group (iii)
applications, she will have the patience to also read the structural theorems.
\ermk

\fac{F5.3} {\bf (The main theorems---sources of more examples).}\ \ 
The following statements are true:
\be
\item 
If $\ct$ has a compact generator $G$, so that $\Hom\big(G,G[n]\big)=0$ for
all $n\geq1$, then $\ct$ is approximable.
\setcounter{enumiv}{\value{enumi}}
\ee
Special cases of (i) include:

The category $\ct=\D(R)$ and the compact generator $G=R$, in other words
we recover Example~\ref{E4.3} as a special case of (i).
More generally: if $R$ is a dga, and $H^n(R)=0$ for all $n>0$, then
the category $\ct=\D(\Mod R)$ with $G=R$ is an example.
Further examples come from topology, for
instance we can let $\ct$ be the homotopy category
of spectra and let $G=\mathbb{S}^0$ be the zero-sphere.

The proof of (i) is basically trivial, there is a
brief discussion in~\cite[Remark~3.3]{Neeman17A}.
\be
\setcounter{enumi}{\value{enumiv}}
\item
Let $X$ be a quasicompact, separated scheme. Then the category
$\Dqc(X)$  is approximable.
\setcounter{enumiv}{\value{enumi}}
\ee
In Remark~\ref{R4.7} we noted that this isn't an easy fact,
it is after all counterintuitive---it says that,
in the category $\Dqc(X)$,
one can pretend to have enough projectives---at least for
some purposes. The proof isn't trivial.

If
X is a separated scheme, of finite type over a noetherian ring,
then the reader can find a proof \cite[Theorem~5.8]{Neeman17}.
It constitutes the main technical lemma of the paper, the rest
amounts to applications.
The generalization to quasicompact, separated schemes is by
a trick which may be found
in \cite[Example~3.6]{Neeman17A}.
\be
\setcounter{enumi}{\value{enumiv}}
\item
Suppose we are given a recollement of triangulated categories
\[\xymatrix@C+20pt{
\car\ar[r] &
\cs\ar@<0.7ex>[l]\ar@<-0.7ex>[l]\ar[r]  &
\ct\ar@<0.7ex>[l]\ar@<-0.7ex>[l]
}\]
with $\car$ and $\ct$ approximable. Assume further that
the category $\cs$ is compactly generated, and any compact object
$H\in\cs$ has the property that $\Hom\big(H,H[i]\big)=0$ for $i\gg0$.
Then the category $\cs$ is also approximable.
\ee
Once again this isn't obvious, it requires proof. The reader can find
it in \cite[Theorem~4.1]{Burke-Neeman-Pauwels18}---it is the
main theorem of the
article.
\efac

So far the majority of the interesting applications has been to algebraic
geometry---it's the example in Fact~\ref{F5.3}(ii) that has proved
useful. But the subject is in its infancy, it is to be hoped that
there will be applications to come, in other contexts.

\fac{F5.5} {\bf (The main theorems---applications).}\ \ 
Assertions  (i) and (ii) below
are~\cite[Theorem~0.5 and Theorem~0.15]{Neeman17},
respectively. Assertion (iii) follows from \cite[Corollary~0.5]{Neeman17A},
while assertion~(iv) follows from \cite[Theorem~0.2]{Neeman18}.
Assertion~(v) is a consequence of \cite[Proposition~0.15]{Neeman18A}, together
with the elaboration and discussion in the couple of paragraphs immediately
following the statement of the Proposition. Anyway: all of (v) may
be found in \cite{Neeman18A}.

In Explanation~\ref{X5.7} the reader is reminded what the various
technical terms in the statements below mean.
\be
\item 
Let $X$ be a quasicompact, separated scheme. The category $\dperf X$
is strongly generated if and only if $X$ has an open cover by affine
schemes $\spec{R_i}$, with each $R_i$ of finite global dimension.
\item
Let $X$ be a separated scheme,
and assume it is noetherian, finite-dimensional, and
that every closed, reduced, irreducible
subscheme of $X$ has a regular alteration.
Then the category 
$\dcoh(X)$ is strongly generated.
\item
Let $X$ be a scheme proper over a noetherian ring $R$. Let $\cy$ be 
the Yoneda map
\[\xymatrix@C+40pt{
\dcoh(X) \ar[r]^-\cy  & \Hom_R^{}\Big(\big[\dperf X\big]\op,\Mod R\Big)
}\]
That is: the map $\cy$ sends the object $B\in\dcoh(X)$ to the functor
$\cy(B)=\Hom(-,B)$, viewed as an $R$--linear homological
functor $\big[\dperf X\big]\op\la\Mod R$.

Then $\cy$
is fully faithful, and the essential image is the set
of finite $R$--linear homological functors $H:\dperf X\op\la \mod R$.
An $R$--linear
homological functor is \emph{finite} if, for all objects
$C\in\dperf X$, the $R$--module $\oplus_nH\big(C[n]\big)$ is finite. 
\item
Suppose $X$ is finite dimensional scheme proper over a noetherian ring $R$,
and assume further that every closed, reduced, irreducible
subscheme of $X$ has a regular alteration.
Let $\wt\cy$ be 
the Yoneda map
\[\xymatrix@C+40pt{
\big[\dperf X\big]\op \ar[r]^-{\wt\cy}  & \Hom_R^{}\Big(\dcoh X,\Mod R\Big)
}\]
That is: the map $\wt\cy$ takes an object
$A\in\dperf X$ to the functor $\wt\cy(A)=\Hom(A,-)$,
viewed as an $R$--linear homological
functor $\dcoh(X)\la\Mod R$.

Then $\wt\cy$ is fully faithful, and the essential image
of $\wt\cy$ are the finite homological functors.
\item
Let $X$ be a noetherian, separated scheme. There is a recipe which
takes the triangulated category $\dperf X$ as input, and out of it constructs
the triangulated category $\dcoh(X)$. And there is a recipe going back:
from the triangulated category $\big[\dcoh(X)\big]\op$
as input the machine spews out
$\big[\dperf X\big]\op$.
\ee
\efac

\xpl{X5.7}
We remind the reader what the terms used in the theorems mean.
Let $\cs$ be a triangulated category, and let $G\in\cs$ be an object.
Then
\be 
\item
$G$ is a \emph{classical generator} if $\cs=\cup_n\genu Gn{-n,n}$.
\item
$G$ is a \emph{strong generator} if there exists an integer
$\ell>0$ with $\cs=\cup_n\genu G\ell{-n,n}$. The category $\cs$ is called
\emph{strongly generated} if it has a strong generator.  
\item
Suppose $X$ is a noetherian scheme,
finite-dimensional, reduced and irreducible. A \emph{regular alteration}
of $X$ is a generically finite, proper, 
surjective morphism $\wt X\la X$ with $\wt X$
regular.
\ee
The non-expert deserves some explanation of (iii):
we all know what a resolution of
singularities is, but the known existence theorems
are too restrictive (for our purposes).
Of course resolutions of singularities are \emph{conjectured} to exist
quite generally, unfortunately what has been proved
so far is limited to equal
characteristic zero, or to schemes of very low dimension.
Regular alterations are less restrictive, and the known existence
theorems are much more
general---see de Jong~\cite{deJong96,deJong97}.

As it turns out, in the proofs of Facts~\ref{F5.5}~(ii) and (iv)
regular alterations suffice.
The non-expert should therefore view the
condition imposed on the noetherian scheme $X$,
in Facts~\ref{F5.5}~(ii) and (iv), as a mild
technical hypothesis.
\expl

\rmk{R5.9}
The reader should note that Fact~\ref{F5.3} asserts that the category
$\Dqc(X)$ is approximable, and now we're telling the reader that the
consequences---Facts~\ref{F5.5}~(i), (ii), (iii), (iv) and (v)---are
all assertions about the categories $\dperf X$ and $\dcoh(X)$.
A technical, formal statement,
about the huge category $\Dqc(X)$, turns out to have
a string of powerful consequences
about the much smaller categories
$\dperf X$ and $\dcoh(X)$, that many people have been
studying for decades.
\ermk

\section{More about the strong generation of $\dperf X$ and $\dcoh(X)$}
\label{S27}

As promised, we will now say a little more about Facts~\ref{F5.5} (i) and (ii).
In this section we will survey what was known before, the basic idea of the
old proofs, and how the proof based on approximability departs from the
older methods.

The non-algebraic-geometers are advised to skip the discussions of the proofs.
The brief summary is that the proofs based on approximability are short, simple,
sweet and work in great generality---the hard work goes into proving that
the category $\Dqc(X)$ is approximable. After that it's all downhill.

Let us begin with Facts~\ref{F5.5}(i), we recall the statement
for the reader's convenience:

\thm{T27.1}
Assume
$X$ is
quasicompact, separated scheme. Then $\dperf X$ is strongly generated if and
only if $X$ may be covered by open affine subsets $\spec{R_i}$, with each
$R_i$ of finite global dimension.
\ethm

\rmk{R27.2}
If $X$ is noetherian and separated, this simplifies to saying that
$\dperf X$ is strongly generated if and only if $X$ is regular and
finite dimensional.
\ermk

\hst{H27.3}
When $X=\spec R$ is affine  Theorem~\ref{T27.1} is old: it
was first proved by Kelly~\cite{Kelly65}, see also
Street~\cite{Street73}.
The result was rediscovered by
Christensen~\cite[Corollary~8.4]{Christensen96}
and later
Rouquier~\cite[Proposition~7.25]{Rouquier08}.

Bondal and Van den Bergh~\cite[Theorem~3.1.4]{BondalvandenBergh04} 
proved the first global
version: if $X$ is a separated scheme, smooth over a field $k$, then 
the category $\dperf X$  is strongly generated. 
The case where $X$ is assumed of finite type over a field and
regular [regularity is weaker than smoothness] follows from
either Rouquier~\cite[Theorem~7.38]{Rouquier08} or
Orlov~\cite[Theorem~3.27]{Orlov16}.

This summarizes the results known before approximability.
Note that, with the exception of Kelly's,
the old results all assumed equal characteristic and that $X$ is noetherian.
By contrast Theorem~\ref{T27.1} works fine in the mixed characteristic,
non-noetherian situation.
\ehst

\pfs{P27.273}
By combining Kelly's old theorem~\cite{Kelly65} with
the main theorem
of Thomason and Trobaugh~\cite{ThomTro},
one easily deduces one of the implications in
Theorem~\ref{T27.1}: if $X$ is quasicompact and separated, and $\spec R$
embeds in $X$ as an open, affine subset, then $R$ must be of finite global
dimension. The reader can
find the argument spelled out in
more detail
in (for example) 
\cite[Remark~0.11]{Neeman17}. 

Now for the tricky direction
of Theorem~\ref{T27.1}, the direction saying that,
if $X$ is quasicompact and separated,
and admits a cover by open affines $\spec{R_i}$ with each $R_i$ of
finite global dimension, then it follows that $\dperf X$ is strongly generated.
As we have already said: the case where $X$ is affine is contained
in Kelly's old theorem.

We remind the reader:
Bondal and Van den Bergh~\cite[Theorem~3.1.4]{BondalvandenBergh04} 
proved the first global
version. They proved that,
if $X$ is a separated scheme, smooth over a field $k$, then 
the category $\dperf X$  is strongly generated. Their proof relies on
the fact that, if $\delta:X\la X\times_kX$ is the diagonal embedding,
then the functor $\R\delta_*$ respects perfect complexes. It is a
characterization of smoothness for 
$\R\delta_*$ to respect perfect complexes---hence 
the argument isn't one that readily lends itself to generalizations.

Nevertheless there were improvements.
The case where $X$ is assumed of finite type over a field and
regular follows from
either Rouquier~\cite[Theorem~7.38]{Rouquier08} or
Orlov~\cite[Theorem~3.27]{Orlov16}. Both proofs still use a diagonal
argument---Rouquier's approach refines Bondal and Van den Bergh's
by stratifying $X$, while the refinement in Orlov's article is
not quite so easy to sum up briefly. It was Orlov's
clever approach to the
problem that inspired the idea of approximability.
\epfs

It remains to give the reader some idea how approximability helps in the
proof of Theorem~\ref{T27.1}. And the main point is that approximability
allows us to reduce the general case to the case of an affine scheme,
where we can use the old theorem of Kelly's. For the reader's convenience
and because the proof is so easy, we prove below the variant of
Kelly's theorem we will actually use.

\thm{T1001.5.7}
Suppose $R$ is an associative ring, and $\D(R)$ its derived category. Let
$n\geq0$ be an integer, and let  
$F\in\D(R)$ be an object such that the projective
dimension of $H^i(F)$ is $\leq n$ for all $i\in\zz$. Then
$F\in\ogen R_{n+1}^{(-\infty,\infty)}$. 
\ethm

Before proving the theorem we remind the reader [who is familiar
  with the calculus of fractions in derived categories]:
any morphism $P\la H^i(E)$ in $\D(R)$, for any
projective $R$--module $P$ and any $E\in\D(R)$, lifts uniquely
to a cochain map
\[\xymatrix@C+10pt{
\cdots \ar[r]& 0  \ar[d]\ar[r]& 0 \ar[d] \ar[r]& P\ar[d]  \ar[r]
&  0 \ar[d]\ar[r]& 0\ar[d] \ar[r] &\cdots\\
\cdots \ar[r]& E^{i-2}  \ar[r]& E^{i-1}  \ar[r]& E^{i}  \ar[r]
&  E^{i+1} \ar[r]& E^{i+2} \ar[r] &\cdots
}\]

\prf
We prove the theorem by induction on $n$. Suppose
first that $n=0$; hence $H^i(F)$ is projective for every
$i\in\zz$. The identity map
$H^i(F)\la H^i(F)$ lifts to a cochain map
\[\xymatrix@C+8pt{
\cdots \ar[r]& 0  \ar[d]\ar[r]& 0 \ar[d] \ar[r]& H^i(F)\ar[d]  \ar[r]
&  0 \ar[d]\ar[r]& 0\ar[d] \ar[r] &\cdots\\
\cdots \ar[r]& F^{i-2}  \ar[r]& F^{i-1}  \ar[r]& F^{i}  \ar[r]
&  F^{i+1} \ar[r]& F^{i+2} \ar[r] &\cdots
}\]
and when we combine, for every $i\in\zz$, we obtain a
cochain map
\[\xymatrix@C-1pt{
\cdots \ar[r]&  H^{-2}(F)  \ar[d]\ar[r]^0&  H^{-1}(F) \ar[d] \ar[r]^0& H^0(F)\ar[d]  \ar[r]^0
&   H^{1}(F) \ar[d]\ar[r]^0&  H^{2}(F)\ar[d] \ar[r] &\cdots\\
\cdots \ar[r]& F^{-2}  \ar[r]& F^{-1}  \ar[r]& F^{0}  \ar[r]
&  F^{1} \ar[r]& F^{2} \ar[r] &\cdots
}\]
This is an isomorphism in cohomology, hence an isomorphism in $\D(R)$.

Now suppose $n\geq0$, and we know the result for every $\ell$ with
$0\leq\ell\leq n$. We wish to show it for $n+1$. Suppose
therefore that we are given an object $F\in\D(R)$ with
$H^i(F)$ of projective dimension $\leq n+1$ for every $i$.
Choose for every $i$ a projective module $P^i$ and a surjection
$P^i\la H^i(F)$. Now form the corresponding cochain map
\[\xymatrix@C+10pt{
\cdots \ar[r]& 0  \ar[d]\ar[r]& 0 \ar[d] \ar[r]& P^{i}\ar[d]  \ar[r]
&  0 \ar[d]\ar[r]& 0\ar[d] \ar[r] &\cdots\\
\cdots \ar[r]& F^{i-2}  \ar[r]& F^{i-1}  \ar[r]& F^{i}  \ar[r]
&  F^{i+1} \ar[r]& F^{i+2} \ar[r] &\cdots
}\]
and combine over $i$ to form
\[\xymatrix@C+13pt{
\cdots \ar[r]& P^{-2}  \ar[d]\ar[r]^0& P^{-1} \ar[d] \ar[r]^0& P^{0}\ar[d]  \ar[r]^0
&  P^{1} \ar[d]\ar[r]^0& P^{2}\ar[d] \ar[r] &\cdots\\
\cdots \ar[r]& F^{-2}  \ar[r]& F^{-1}  \ar[r]& F^{0}  \ar[r]
&  F^{1} \ar[r]& F^{2} \ar[r] &\cdots
}\]
giving a map $P\la F$, which we complete to a triangle $P\la F\la Q$.
Clearly $P\in\ogen R_{1}^{(-\infty,\infty)}$, and the long
exact sequence in cohomology gives that $H^i(Q)$ is of projective
dimension $\leq n$. Hence $F$ belongs to
$\ogen R_{1}^{(-\infty,\infty)}*\ogen R_{n+1}^{(-\infty,\infty)}
\subset\ogen R_{n+2}^{(-\infty,\infty)}$.
\eprf

At this point the non-algebraic-geometer (who hasn't yet
done so) is advised to skip ahead to
Theorem~\ref{T27.17}. What will come between now and then is largely
aimed to show that the approximability of $\Dqc(X)$ makes the
reduction to Kelly's old theorem straightforward and easy---hopefully
the sketch we give will make this transparent to the experts, but
for non-algebraic-geometers it might be mystifying. Anyway:
the reduction depends on the following little lemma---the reader
should note the way
approximability enters the proof of the lemma, this is
the only point where approximability will be used.

\lem{L27.5}
Let $X$ be a quasicompact, separated scheme, let $G\in\Dqc(X)$
be a compact generator, 
and let $u:U\la X$
be an open immersion with $U$ quasicompact. Then the object 
$\R u_*\co_U^{}\in\Dqc(X)$ belongs to $\ogenu Gn{-n,n}$ for
some integer $n>0$.
\elem

\prf
It is relatively easy to show that, for some sufficiently large integer
$\ell>0$, we have $\Hom\big(\R u_*\co_U^{}\,\,,\,\,\Dqc(X)^{\leq-\ell}\big)=0$.
By the approximability\footnote{This isn't immediate from the definition
of approximability, but follows from the structural theorems.
We are using the fact that
$\R u_*\co_U^{}\in\Dqc(X)^{\leq m}$ for some $m>0$, coupled with the
fact that one can approximate objects in $\Dqc(X)^{\leq0}$ to arbitrary
order, not just to order zero as given in the definition.
See Sketch~\ref{S29.29}(i) for more detail.} of $\Dqc(X)$ we
may choose an integer $n$ and a triangle 
$E\la\R u_*\co_U^{}\la D$ with $D\in\Dqc(X)^{\leq-\ell}$ and
$E\in\ogenu Gn{-n,n}$.

But the map $\R u_*\co_U^{}\la D$ must vanish by the choice of $\ell$,
making $\R u_*\co_U^{}$ a direct summand of the object 
$E\in\ogenu Gn{-n,n}$.
\eprf

\skt{S27.7}
We should indicate how Theorem~\ref{T27.1} follows from
the combination of Lemma~\ref{L27.5}
and Kelly's old theorem.
Let $X$ and the open affine cover by $U_i=\spec{R_i}$ be
as in the hypotheses of
Theorem~\ref{T27.1}. Because $X$ is quasicompact we may, possibly after
passing to a subcover, assume that our cover is
finite; write the cover as $\{U_i,\,1\leq i\leq r\}$.

Now choose
a compact generator $G\in\Dqc(X)$. The Lemma allows us
to choose, for each open subset $U_i$, an integer $n_i$ so that 
$\R u_{i*}\co_{U_i}^{}\in\ogenu G{n_i}{-n_i,n_i}$. Let $n$ be the
maximum of the finitely many $n_i$; then
$\R u_{i*}\co_{U_i}^{}\in\ogenu G{n}{-n,n}\subset\ogen G_n^{(-\infty,\infty)}$
for every $i$ in the finite
set.

Next, for each $i$ we know that
$U_i=\spec{R_i}$ with $R_i$ is of finite global dimension, and
Theorem~\ref{T1001.5.7} tells us
that we may choose an integer $\ell>0$ so that,
for every one of the finitely many $i$ with $1\leq i\leq r$, we
have
$\Dqc(U_i)=\ogen {\co_{U_i}^{}}_\ell^{(-\infty,\infty)}$. It follows
that
\[
\R u_{i*}\Dqc(U_i)
\quad=\quad\R u_{i*}\Big[\ogen {\co_{U_i}^{}}_\ell^{(-\infty,\infty)}\Big]
\quad\subset\quad\ogen {\R u_{i*}\co_{U_i}^{}}_\ell^{(-\infty,\infty)}
\quad\subset\quad\ogen G_{\ell n}^{(-\infty,\infty)}
\]
Let $\cv=\add\big[\cup_{i=1}^{r}\R u_{i*}\Dqc(U_i)\big]$,
with the notation as in Reminder~\ref{R3.1}(ii). By
the displayed inclusion above
$\big[\cup_{i=1}^{r}\R u_{i*}\Dqc(U_i)\big]\subset
\ogen G_{\ell n}^{(-\infty,\infty)}$,
and as $\ogen G_{\ell n}^{(-\infty,\infty)}$
is closed under (finite)
coproducts it
follows that $\cv\subset\ogen G_{\ell n}^{(-\infty,\infty)}$.

It's an exercise to show that
\[
\Dqc(X)\eq \underbrace{\cv*\cv*\cdots*\cv}_{r\text{ copies}}
\]
with the notation as in Reminder~\ref{R3.1}(i). Hence
$\Dqc(X)=\ogen G_{\ell nr}^{(-\infty,\infty)}$.
We have proved a statement about $\Dqc(X)$, and
in Notation~\ref{N1.107589} we learned that  $\dperf X$ is equal to
the subcategory of compact objects in $\Dqc(X)$. Standard compactness
arguments tell us that from the equality
$\Dqc(X)=\ogen G_{\ell nr}^{(-\infty,\infty)}$
we can formally deduce the equality
$\dperf X=\cup_{m>0}^{}
\genu G{\ell nr}{-m,m}$.\hfill{$\Box$}
\eskt

We want to highlight the power of approximability.
Sketch~\ref{S27.7} was meant to show the expert that
Theorem~\ref{T27.1} is easy to deduce by combining
Kelly's old theorem with Lemma~\ref{L27.5}, and 
the proof of Lemma~\ref{L27.5} displays
how the lemma follows immediately from the fact that
$\Dqc(X)$ is approximable.

While we're into exhibiting the power of
approximability, 
let us mention another corollary
of Lemma~\ref{L27.5}---and therefore another easy
consequence of approximability.

\thm{T27.9}
Suppose $f:X\la Y$ is a separated morphism of
quasicompact, quasiseparated schemes. If
$\R f_*:\Dqc(X)\la\Dqc(Y)$ takes perfect complexes to complexes
of bounded--below Tor-amplitude then $f$ must be of finite Tor-dimension.
\ethm

\rmd{R27.11}
We owe the reader a glossary of the technical terms in the statement of
Theorem~\ref{T27.9}.
\be 
\item  
Given a morphism of schemes $f:X\la Y$, for any $x\in X$ there is
an induced ring homomorphism $\co_{Y,f(x)}^{}\la\co_{X,x}^{}$ of the 
stalks. The map $f$ is \emph{of finite Tor-dimension at $x$} if 
$\co_{X,x}^{}$ has a finite flat resolution over $\co_{Y,f(x)}^{}$. 
\item
The map $f$ is \emph{of finite Tor-dimension} if it is of finite Tor-dimension
at every $x\in X$.
\item
The complex $C\in\Dqc(Y)$ is \emph{of bounded-below Tor-amplitude} if,
for every open immersion $u:U\la Y$ with $U=\spec R$ affine,
the complex $u^*C\in\Dqc(U)\cong\D(R)$ is
isomorphic to a bounded-below $K$--flat complex.
\ee 
\ermd

\hst{H27.15}
We should tell the reader what was known in the
direction of Theorem~\ref{T27.9}.
If the schemes $X$ and $Y$ are noetherian and $f:X\la Y$ is of finite
type, then the converse of Theorem~\ref{T27.9} is known and old---the
reader may find it in Illusie~\cite[Corollaire~4.8.1]{Illusie71B}.
The direction proved in Theorem~\ref{T27.9} was open for a long time,
the first progress was in \cite{Lipman-Neeman07}.
But the statement in
\cite{Lipman-Neeman07} is much narrower than  Theorem~\ref{T27.9}, it
is confined to the situation where $f$ is proper.
\ehst

\skt{S27.13} We should give the reader some idea why Theorem~\ref{T27.9}
follows easily from Lemma~\ref{L27.5}---this discussion is
for algebraic geometers, the non-specialists are advised to skip
ahead to Theorem~\ref{T27.17}.

It's obviously local in $Y$ to determine if $f$ is of finite Tor-dimension.
Using the main theorem of Thomason and Trobaugh~\cite{ThomTro},
it's also local in $Y$ to determine whether $\R f_*$ takes
perfect complexes to complexes of bounded-below Tor-amplitude.
Hence we may assume $Y$ is affine,
therefore separated. As $f$ is separated we deduce that $X$ must be separated.  

We are given that $\R f_*$ takes perfect complexes to
complexes of bounded-below Tor-amplitude, and wish to show
that $f$ is of finite Tor-dimension. Being of finite Tor-dimension is
local in $X$; it suffices
to show that, for each of open immersion $u:U\la X$ with $U$ affine,
the composite
$U\stackrel u\la X\stackrel f\la Y$
is of finite Tor-dimension. By Lemma~\ref{L27.5} there exists a perfect
complex $G\in\Dqc(X)$ and an integer $n>0$ with
$\R u_*\co_U^{}\in\ogenu Gn{-n,n}$. Therefore
\[\begin{array}{ccccc}
(fu)_*\co_U^{}&\quad\cong\quad&
\R(fu)_*\co_U^{}&\quad\cong\quad&
\R f_*\R u_*\co_U^{} \\*[5pt]
&\quad\in\quad&
\R f_*\Big[\ogenu Gn{-n,n}\Big]&\quad\subset\quad&
\ogenu {\R f_*G}n{-n,n}
\end{array}\]
But $\R f_*G$ is of bounded below Tor-amplitude by hypothesis, and
in forming $\ogenu {\R f_*G}n{-n,n}$
we only allow
$\R f_*G[i]$ with $-n\leq i\leq n$, coproducts, extensions and direct summands.
Hence 
the objects of $\ogenu {\R f_*G}n{-n,n}$ have Tor-amplitude uniformly
bounded below.
\eskt

It's time to turn our attention to Fact~\ref{F5.5}(ii), we remind the
reader of the statement:

\thm{T27.17}
Let $X$ be a separated, noetherian, finite-dimensional scheme,
and assume that every closed, reduced, irreducible
subscheme of $X$ has a regular alteration.
Then the category 
$\dcoh(X)$ is strongly generated.
\ethm

\hst{R27.21}
We should tell the reader what was known in
the direction of Theorem~\ref{T27.17}.
We have already alluded to the fact that,
when $X$ is regular and finite-dimensional, 
the inclusion $\dperf X\la\dcoh(X)$ is an equivalence
and Theorem~\ref{T27.1} tells us that the equivalent categories
$\dperf X\cong\dcoh(X)$ are strongly generated.
Using a stratification, of a possibly singular $X$,
Rouquier~\cite[Theorem~7.38]{Rouquier08}
built and substantially extended on the argument in Bondal and
Van den Bergh~\cite[Theorem~3.1.4]{BondalvandenBergh04}
to show that $\dcoh(X)$ is strongly generated whenever $X$ is a 
separated scheme 
of finite type over a perfect field $k$. The preprint by
Keller and 
Van den Bergh~\cite[Proposition~5.1.2]{Keller-Vandenbergh08}
generalized to separated schemes of finite type over arbitrary fields,
but this Proposition disappeared in the passage to the published 
version~\cite{Keller-Vandenbergh11}.
The reader might also wish to look at Lunts~\cite[Theorem~6.3]{Lunts10}
for a different approach to the proof, but still using stratifications.
If we specialize the result of Rouquier, extended by Keller and 
Van den Bergh, to the case where $X=\spec R$ is an affine scheme,
we learn that $\D^b(\mmod R)$ is strongly generated whenever $R$ is of finite
type over a field $k$.

Note that, while Theorem~\ref{T27.1} is easy and classical in the
case where $X$ is affine, Theorem~\ref{T27.17} is \emph{neither easy nor
classical for affine $X$.}
In recent years there has been interest 
among commutative algebraists in understanding this better: 
the reader is referred
to 
Aihara and Takahashi~\cite{Aihara-Takahashi11},
Bahlekeh, 
Hakimian, Salarian
and Takahashi~\cite{Bahlekeh-Hakimian-Salarian-Takahashi15} and
Iyengar and Takahashi~\cite{Iyengar-Takahashi15} 
for a sample of the literature. There is also a connection
with the concept of the radius of the (abelian) category of
modules over $R$; see Dao and Takahashi~\cite{Dao-Takahashi14,Dao-Takahashi15}
and Iyengar and Takahashi~\cite{Iyengar-Takahashi15}.
The union of the known results
seems to be that $\D^b(\mmod R)$ is
strongly generated if $R$ is an
equicharacteristic  excellent
local ring, or essentially
of finite type over a field---see~\cite[Corollary~7.2]{Iyengar-Takahashi15}.
In~\cite[Remark~7.3]{Iyengar-Takahashi15} it is observed
that there are examples of commutative, noetherian 
rings for which $\D^b(\mmod R)$ is \emph{not} strongly generated.

The structure of the proof
of Theorem~\ref{T27.17}
(see Sketch~\ref{S27.19}) is
that one passes to regular alterations of $X$ and its closed
subschemes. Assuming $X$ affine is no help with
the approximability proof of
Theorem~\ref{T27.17}---when $X$ is affine
and singular
we end up proving
a result in commutative algebra, but the technique of the proof
passes through non-affine schemes. 

Unlike all the pre-approximability results, except Kelly's,
Theorems~\ref{T27.1} and \ref{T27.17} do
not assume equal characteristic.
\ehst

\skt{S27.19}
We should tell the reader a little about the proof. But first we should
make it clear that Theorem~\ref{T27.17} will not be proved using
approximability directly, instead we will prove it as a corollary
of Theorem~\ref{T27.1}, which followed from approximability.
Precisely: Theorem~\ref{T27.1} and Theorem~\ref{T27.17} are identical
when $X$ is a finite-dimensional, regular, noetherian, separated scheme. And 
the idea is to reduce to this case.

Resolutions of singularities might look tempting, but in mixed characteristic
they are known to exist only in low dimension. So the key is that
we can get by with regular alterations---the hypotheses of
the theorem say that they exist for every closed subvariety of $X$,
and it turns
out that Theorem~\ref{T27.17} can be deduced from this using induction on the
dimension of $X$ and two old theorems of Thomason's.

This survey has been stressing that the hard work goes into proving
approximability, the consequences are all easy corollaries. Theorem~\ref{T27.17}
must count as an exception, the argument is tricky. It might be relevant
to note that in this field---noncommutative algebraic geometry---there
are quite a number of theorems that are known in characteristic zero with
proofs that rely on resolutions of singularities, and
conjectured in positive characteristic. I wasn't the first to come up with
the idea of trying to use de Jong's theorem, in other words trying to
prove these conjectures using regular alterations. So far
Theorem~\ref{T27.17} is the only success story. It isn't regular alterations
alone that do the trick, it's the combination of regular alterations and
support theory---in this case support theory manifests itself as the two
old theorems of Thomason's.
\eskt

\plm{P27.995}
There is a non-commutative version---Kelly's old theorem doesn't
assume commutativity. This raises the obvious question: to what extent
do the more recent
theorems extend beyond commutative algebraic geometry?

Perhaps we should explain, and for simplicity let us stick to the
case where $X=\spec R$ is affine. As we have presented the
theory, up to now, we have
implicitly been assuming that the ring $R$ is commutative. But what Kelly
proved doesn't depend on commutativity---the reader can see
this for herself, just look at the proof
of Theorem~\ref{T1001.5.7}. Let $R$ be any associative ring
and let $\D^b(\proj R)$ be the derived category of bounded complexes of
finitely generated, projective $R$--modules. Kelly's 1965 theorem says
that $D^b(\proj R)$ has a strong generator if and only if $R$ is of finite
global dimension.

All the later theorems listed above, including the recent ones
whose proof 
relies on approximability, assume commutativity. In particular:
assume $R$ is a commutative, noetherian ring,
of finite type
over an excellent ring of dimension $\leq2$.
Theorem~\ref{T27.17}, in the special case where $X=\spec R$,
tells us that the category
$\dcoh(X)\cong\D^b(\mod R)$ is strongly generated. The category
$\D^b(\mod R)$ has for objects the bounded complexes of finite $R$--modules.

Is the commutativity hypothesis necessary in the above? Is there some
large class of noncommutative, noetherian rings for which $\D^b(\mod R)$ is
strongly generated? The proof in the commutative case, which goes by way
of the regular alterations of de Jong, doesn't seem capable of a noncommutative
extension.
\eplm

\rmk{R27.101978}
Recall: a strong generator in $\cs$ is an object $G\in\cs$ such that,
for some integer $\ell>0$, we have
$\cs=\cup_n\genu G\ell{-n,n}$. One can ask for estimates on
$\ell$.
This leads to the definitions

\be
\item
Given objects $G,F\in\cs$, the \emph{$G$--level of $F$} is the smallest
integer $\ell$ such that $F\in \cup_n\genu G{\ell+1}{-n,n}$. This
notion is due to Avramov, Buchweitz and 
Iyengar~\cite{Avramov-Buchweitz-Iyengar07}.
\item
Let $G$ be an object of $\cs$. The \emph{generation time} of $G$ is the
smallest $\ell$ for which $\cs= \cup_n\genu G{\ell+1}{-n,n}$.
The set of all possible generation times, taken over all strong
generators $G\in\cs$, is known as the \emph{Orlov spectrum} of $\cs$.
These notions first appeared in Orlov~\cite{Orlov09}.
\item
The \emph{Rouquier dimension} of $\cs$ is the smallest integer $\ell$ such 
that there exists a $G$ with $\cs= \cup_n\genu G{\ell+1}{-n,n}$.
This integer first appeared in Rouquier~\cite{Rouquier08}---Rouquier's
name for this number was just plain ``dimension''.
\setcounter{enumiv}{\value{enumi}}
\ee
There are several conjectures, and many papers estimating these
numbers---almost all in the equal charateristic case, after all until
recently there was no existence theorem of strong generators
in mixed characteristic.
One can ask if the theorems surveyed in this article give good
bounds in mixed characteristic---and the short answer is No. In more detail:
\be
\setcounter{enumi}{\value{enumiv}}
\item
If we assume that $X$ is regular and
quasiprojective, then the proof
of Application~1(i) is effective. It gives an explicit upper bound on
the
Rouquier dimension of $\dperf X=\dcoh(X)$. But the bound is dreadful.
\item
If we drop the quasiprojectivity hypothesis, and/or if we allow
singularities, then the proof becomes ineffective. It proves the
existence of an integer $\ell>0$ and a generator $G$ with 
$\dcoh(X)= \cup_n\genu G\ell{-n,n}$, but there is no estimate on $\ell$.
\ee
\ermk

\ela{E27.9098325}
The following is for the benefit of the readers who would like
Remark~\ref{R27.101978}
spelt out a little more.

In~\cite[Section~4]{Neeman17} the reader can find the argument that
proves the approximability of $\Dqc(X)$ when $X$ is
quasiprojective---and
in  \cite[Proposition~4.4]{Neeman17} it's made clear that the
estimates are explicit. And, assuming $X$ is not only
quasiprojective but also regular,  Sketch~\ref{S27.7} shows us how to pass
from the estimates given by approximability to explicit estimates on
the $\ell$ for which $\dperf X=\cup_n\genu G{\ell+1}{-n,n}$. Still
assuming that $X$ is quasiprojective and regular, a careful reading of
the
proof of \cite[Proposition~4.4]{Neeman17} will show us that these
crude
estimates can easily be improved. But our point for now is 
that the proof is effective, it gives bounds.

The general proof of approximability, for quasicompact, separated $X$,
is by reduction to the quasiprojective case. This reduction goes by
\be
\item
We first use noetherian approximation to reduce to the case where $X$
is of finite type over $\zz$.
\item
Next we use induction on the dimension coupled with Chow's Lemma to 
reduce to the quasiprojective case.  
\setcounter{enumiv}{\value{enumi}}
\ee
The way to use 
noetherian approximation is straightforward enough---we don't
go into detail, leaving this to the reader's imagination. In
principle this part could be made effective, but not in a 
way that yields good bounds. Still: so far there are bounds of some sort. 

But the true subtlety arises in (ii), with the way we use Chow's Lemma. 
Chow's Lemma produces
for us a birational, projective morphism $\pi:\wt X\la X$ with $\wt X$
quasiprojective.
And two remarks are in order
\be
\setcounter{enumi}{\value{enumiv}}
\item
Chow's Lemma is where we have to assume $X$ separated. This is the
point of the proof where it doesn't suffice for $X$ to be quasiseparated.
\item
Since $\wt X$ is quasiprojective the category $\Dqc(\wt X)$ is
approximable. To deduce from this useful information about $\Dqc(X)$
we
need to take a compact generator $\wt G\in\Dqc(\wt X)$, and approximate
$\R\pi_*\wt G$ using a compact generator $G\in\Dqc(X)$.
\ee
There is a way to achieve (iv), it relies on 
\cite[Theorem~4.1]{Lipman-Neeman07}. But the proof of 
\cite[Theorem~4.1]{Lipman-Neeman07} is homotopy-theoretic, it hinges on
the fact that, for a compact object $G$, the functor $\Hom(G,-)$ commutes
with homotopy colimits. Ignoring the technicality that we work with 
homotopy colimits and not ordinary colimits, the point is the following.
Any map from a compact object $G$ to an object $\colim\, T_i$ factors through
some $T_i$, \emph{but we have no control over 
how large $i\in\nn$ will have to be.} Thus any argument which has, embedded in
it, the appeal to such homotopy colimit arguments, is inherently and
hopelessly ineffective.
This explains why we lose control when $X$ isn't assumed quasiprojective.

If we allow $X$ to become singular, then  Sketch~\ref{S27.19} hints
how to reduce to the regular case using regular alterations. 
As was said in Sketch~\ref{S27.19}:
the proof appeals to 
two old theorems of Thomason's, both of which depend 
on
homotopy-colimit arguments. Hence this passage is also ineffective 
beyond salvation.
\eela

\section{More about finite $R$--linear functors
  $H:\big[\dperf X\big]\op\la\Mod R$ and $\wt H:\dcoh(X)\la\Mod R$}
\label{S29}

It's time to expand on Facts~\ref{F5.5}~(iii) and (iv). We begin by recalling 
the
statements.

\thm{T29.1}
Let $X$ be a finite-dimensional
scheme proper over a noetherian ring $R$. Let $\cy$
the Yoneda map
\[\xymatrix@C+40pt{
\dcoh(X) \ar[r]^-\cy  & \Hom_R^{}\Big(\dperf X\op,\Mod R\Big)
}\]
taking $B\in\dcoh(X)$ to the functor $\Hom(-,B)$, and
let $\wt\cy$ be 
the Yoneda map
\[\xymatrix@C+40pt{
\big[\dperf X\big]\op \ar[r]^-{\wt\cy}  & \Hom_R^{}\Big(\dcoh X,\Mod R\Big)
}\]
taking $A\in\dperf X$ to the functor $\Hom(A,-)$.
Assuming every closed subvariety of $X$ admits a regular alteration
both functors
are fully faithful, and in each case the essential image is the set
of finite $R$--linear homological functors.
Recall: an $R$--linear
homological functor $H:\cs\la\Mod R$ is \emph{finite} if, for all objects
$C\in\cs$, the $R$--module $\oplus_nH\big(C[n]\big)$ is finite.

For the functor $\cy$ the assertion is true even without the hypotheses
of finite-dimensionality and
the existence of regular alterations.
\ethm

\hst{R29.3}
We remind the reader what was known before.
\be 
\item
If $X$ is proper over $R$, if $A\in\dperf X$ and if $B\in\dcoh(X)$, then
\[
\Hom(A[i],B)\quad\cong\quad H^{-i}(A^\vee\oo B)
\]
is a finite $R$--module for every $i$ and vanishes outside a bounded range.
This much was proved by Grothendieck~\cite[Th\'eor\`eme~3.2.1]{Grothendieck61}.

Translating to the language of Theorem~\ref{T29.1}: given objects
$A\in\dperf X$ and $B\in\dcoh(X)$, then
$\cy(B)=\Hom(-,B)$ is a finite
homological functor on $\big[\dperf X\big]\op$, while $\wt\cy(A)=\Hom(A,-)$ is a
finite homological functor on $\dcoh(X)$.
This much has been known since 1961.
\item
\emph{As long as $R$
is
a field,} Bondal and Van den Bergh~\cite[Theorem~A.1]{BondalvandenBergh04}
proved that 
every finite homological functor on $\big[\dperf X\big]\op$ is
$\Hom(-,B)$ for some $B\in\dcoh(X)$. In the language of
Theorem~\ref{T29.1}: they proved that the essential image of
$\cy$ consists of the finite homological functors.
\item
\emph{Still assuming $R$
is
a field,} the assertion of
Theorem~\ref{T29.1} about the functor $\wt\cy$
can be found in
Rouquier~\cite[Corollary~7.51(ii)]{Rouquier08}---although
the author of the present article doesn't follow 
the argument in~\cite{Rouquier08} that briefly
outlines how a proof might go, it's too skimpy.
\ee
If $R$ is a field Theorem~\ref{T29.1} improves on what was known
about the functor $\cy$
by showing that it's fully faithful. And for $R$ more general
Theorem~\ref{T29.1} is new, for both the functor $\cy$ and the functor
$\wt\cy$.
\ehst

And now the time has come to tell the reader something about the proof
of Theorem~\ref{T29.1}. It turns out that the theorem is an immediate
corollary of a far more general fact, and the discussion
of this result brings us naturally to the structure that all approximable
categories share. Let us begin in even greater generality, not assuming
all the hypotheses of approximability.

\dfn{D29.5}
Let $\ct$ be a triangulated category, and let
$\big(\ct^{\leq0}_1,\ct^{\geq0}_1\big)$ and $\big(\ct^{\leq0}_2,\ct^{\geq0}_2\big)$
be two {\it t}--structures on $\ct$. We declare them \emph{equivalent}
if there exists an integer $A>0$ with
$\ct^{\leq-A}_1\subset\ct^{\leq0}_2\subset\ct^{\leq A}_1$.
\edfn

The definition agrees with the intuition of the Introduction: each \tstr\
defines a kind of (directed)
metric, and we'd like declare {\it t}--structures
equivalent whenever
they induce equivalent metrics. And now we recall

\rmk{R29.7}
Let $\ct$ be a triangulated category with coproducts, and let $G\in\ct$ be a
compact generator.
From Alonso, Jerem{\'{\i}}as and
Souto~\cite[Theorem~A.1]{Alonso-Jeremias-Souto03} we
learn that $\ct$ has a unique \tstr\ $\big(\ct_G^{\leq0},\ct_G^{\geq0}\big)$
\emph{generated by $G$.}

It is not difficult to show that, if $G$ and $H$ are two compact generators
for $\ct$, then the {\it t}--structures $\big(\ct_G^{\leq0},\ct_G^{\geq0}\big)$
and $\big(\ct_H^{\leq0},\ct_H^{\geq0}\big)$ are equivalent. Thus
\emph{up to equivalence} there is a preferred \tstr\ on $\ct$, namely
$\big(\ct_G^{\leq0},\ct_G^{\geq0}\big)$ where $G$ is a compact generator.
We say that a \tstr\
$\big(\ct^{\leq0},\ct^{\geq0}\big)$ is in the \emph{preferred equivalence
class} if it is equivalent to $\big(\ct_G^{\leq0},\ct_G^{\geq0}\big)$ for
some compact generator $G$, hence for every compact generator.
\ermk

\dis{D29.9}
Given a \tstr\ $\big(\ct^{\leq0},\ct^{\geq0}\big)$ it is customary to define
the 
categories
\[
\ct^-=\cup_n\ct^{\leq n}\,,\qquad \ct^+=\cup_n\ct^{\geq -n}\,,
\qquad \ct^b=\ct^-\cap\ct^+
\]
as in Notation~\ref{N1.23}. It's obvious from Definition~\ref{D29.5}
that equivalent {\it t}--structures yield identical
$\ct^-$, $\ct^+$ and $\ct^b$.

Now assume we are in the situation of Remark~\ref{R29.7}, that is $\ct$
has coproducts and there exists a single compact generator $G$.
Then there is a preferred equivalence class of {\it t}--structures
and, correspondingly, preferred $\ct^-$, $\ct^+$ and $\ct^b$. These
are intrinsic, they're independent of any choice. In the remainder
of the article we only consider the ``preferred''
$\ct^-$, $\ct^+$ and $\ct^b$.
\edis

Slightly more sophisticated is the category $\ct^-_c$ below.

\dfn{D29.11}
Let $\ct$ be a triangulated category with coproducts, and assume it has
a compact generator $G$. Choose a \tstr\ $\big(\ct^{\leq0},\ct^{\geq0}\big)$
in the preferred equivalence class.
The full subcategory $\ct^-_c$ is defined by
\[
\ct^-_c\eq\left\{F\in\ct\left|\begin{array}{c}
\text{For all integers }n>0\text{ there exists a triangle}\\
E\la F\la D\la E[1]\\
\text{with $E$ compact and }D\in\ct^{\leq-n-1}
\end{array}\right.\right\}
\]
We furthermore define $\ct^b_c=\ct^b\cap\ct^-_c$.
\edfn

\rmk{R29.13}
Intuitively the category $\ct^-_c$ is the
closure, with respect to the metric induced
by the \tstr\ $\big(\ct^{\leq0},\ct^{\geq0}\big)$,
of the subcategory $\ct^c$ of all compact objects in $\ct$.
It's obvious that the category $\ct^-_c$ is intrinsic, after
all equivalent metrics will lead to the same closure. And
as $\ct^-_c$ and $\ct^b$ are both intrinsic, so is their
intersection $\ct^b_c$.
\ermk

We have defined all this intrinsic structure, assuming only that
$\ct$ is a triangulated category with coproducts and with a single
compact generator. In this generality we know that the
subcategories $\ct^-$, $\ct^+$ and $\ct^b$ are thick. For the
subcategories $\ct^-_c$ and $\ct^b_c$ one proves

\pro{P29.15}
If $\ct$ has a compact generator $G$, such that $\Hom\big(G,G[n]\big)=0$
for $n\gg0$, then the subcategories $\ct^-_c$ and $\ct^b_c$ are
thick.
\epro

\rmk{R29.16}
If $\ct$ is approximable then, by Definition~\ref{D4.1}(i), there is
an integer $A>0$, a compact generator $G\in\ct$
and a \tstr\ $\big(\ct^{\leq0},\ct^{\geq0}\big)$, so that
$\Hom\big(G[-A],\ct^{\leq0}\big)=0$
and $G[A]\in\ct^{\leq0}$. Hence $G[n]\in\ct^{\leq0}$ for all $n\geq A$
and
$\Hom\big(G,G[n]\big)=0$ for all $n\geq2A$; Proposition~\ref{P29.15}
therefore tells us that the subcategories $\ct^-_c$ and $\ct^b_c$ are
thick whenever $\ct$ is approximable.
\ermk

Of course it would be nice to be able to work out examples: what does
all of this
intrinsic structure come down
to in special cases? This is where approximability
helps. We first note

\pro{P29.17}
Assume the category $\ct$ is approximable; see Definition~\ref{D4.1}.
We recall part of the definition: the category $\ct$ is approximable
if it has a compact generator $G$, a \tstr\ $\big(\ct^{\leq0},\ct^{\geq0}\big)$
and an integer $A>0$
satisfying some properties, see Definition~\ref{D4.1}~(i) and (ii) for
the properties.

Then any \tstr, which comes as part of a triad
satisfying the properties of Definition~\ref{D4.1}~(i) and (ii),
must be in the preferred equivalence
class. Furthermore: for any compact generator $G'$ and any
\tstr\ $\big(\ct^{\leq0},\ct^{\geq0}\big)$ in the preferred equivalence class,
there is an integer $A'>0$ so that the properties of
Definition~\ref{D4.1}~(i) and (ii) hold.
\epro

In practice this means that, in proving that $\ct$ is approximable, we must
produce at least one useful \tstr\ that we know belongs to the preferred
equivalence class. After all: this \tstr\
must be manageable enough
to lend itself to a proof that the conditions in
Definition~\ref{D4.1}~(i) and (ii) hold.
Note that the proof of
Alonso, Jerem{\'{\i}}as and
Souto~\cite[Theorem~A.1]{Alonso-Jeremias-Souto03}
yields a 
\tstr\ $\big(\ct_G^{\leq0},\ct_G^{\geq0}\big)$ in the
preferred equivalence class, but the construction is a little
opaque---it shows existence and uniqueness,
but usually doesn't give us much of
a handle on $\big(\ct_G^{\leq0},\ct_G^{\geq0}\big)$. So while we know
that {\it t}--structures in the preferred equivalence class exist, this
needn't be especially useful in working with them.

Let $X$ be a quasicompact, separated scheme. We have told the reader that
\cite[Theorem~5.8]{Neeman17} combined
with \cite[Example~3.6]{Neeman17A} prove that $\ct=\Dqc(X)$ is approximable;
the \tstr\ used in the proof happens to be the standard \tstr\
of Example~\ref{E1.17}~(iii) and (iv). We remind
the reader: the standard \tstr\ on $\Dqc(X)$ has
\begin{eqnarray*}
\Dqc(X)^{\leq0} &=& \{F\in\Dqc(X)\mid \ch^i(F)=0\text{ for all }i>0\}\\
\Dqc(X)^{\geq0} &=& \{F\in\Dqc(X)\mid \ch^i(F)=0\text{ for all }i<0\}
\end{eqnarray*}
where $\ch^i$ is the functor taking a cochain complex to its
$i\mth$ cohomology sheaf.
Proposition~\ref{P29.17} now informs us that the standard \tstr\ must
belong to the preferred equivalence class. Hence the
categories $\ct^-$, $\ct^+$ and $\ct^b$ are the usual: we have
$\ct^-=\Dqcmi(X)$, $\ct^+=\Dqcpl(X)$ and $\ct^b=\Dqcb(X)$. 
The subcategories $\Dqcmi(X)$, $\Dqcpl(X)$ and $\Dqcb(X)$ of $\Dqc(X)$
are traditionally defined to be
\begin{eqnarray*}
\Dqcmi(X)&=& \{F\in\Dqc(X)\mid \ch^i(F)=0 \text{ for all }i\gg0\} \\
\Dqcpl(X)&=& \{F\in\Dqc(X)\mid \ch^i(F)=0 \text{ for all }i\ll0\} \\
\Dqcb(X)&=&\Dqcmi(X)\cap\Dqcpl(X)
\end{eqnarray*}

Next we ask ourselves: what about $\ct^-_c$ and $\ct^b_c$? We begin with the
affine case.

\exe{E29.19}
Let $R$ be a ring. Prove that, in the category $\ct=\D(R)$, the subcategory
$\ct^-_c$ agrees with the $\D^-(\proj R)$ of Remark~\ref{R3.13}.
\eexe

\obs{O29.21}
Now assume $R$ is a commutative ring, and let $X=\spec R$. Then
\cite[Theorem~5.1]{Bokstedt-Neeman93} tells us that the
natural functor
$\D(R)\la\Dqc(X)$ is an equivalence of categories. Putting
$\ct=\Dqc(X)\cong\D(R)$, we learn from Exercise~\ref{E29.19}
what the category $\ct^-_c$ is.

Now let $X$ be any quasicompact, separated scheme. If $u:U\la X$ is an
open immersion, then the functor $u^*:\Dqc(X)\la\Dqc(U)$ respects the standard
\tstr\ and sends compact objects in $\Dqc(X)$ to compact objects
in $\Dqc(U)$. Hence $u^*\Dqc(X)^-_c\subset\Dqc(U)^-_c$. Thus
every object in $\Dqc(X)^-_c$ must be ``locally in $\D^-(\proj R)$'',
meaning that for every open immersion $u:\spec R\la X$ we must have
that $u^*\Dqc(X)^-_c\subset\D^-(\proj R)$. The objects ``locally in
$\D^-(\proj R)$'' were first studied by
Illusie~\cite{Illusie71B,Illusie71A} in SGA6. They have a name,
they are the pseudocoherent complexes.

The next result is not so obvious. In \cite[Theorem~4.1]{Lipman-Neeman07}
the reader can find a proof that
\eobs

\pro{P29.23} Let $X$ be a quasicompact, separated scheme.
Then every pseudocoherent complex belongs to $\Dqc(X)^-_c$. Coupled
with Observation~\ref{O29.21} this teaches us that the objects
of $\Dqc(X)^-_c$ are precisely the pseudocoherent complexes.
\epro

\rmk{R29.25}
From now on we will assume the scheme $X$ noetherian and separated.
In this case pseudocoherence simplifies: we have $\Dqc(X)^-_c=\dmcoh(X)$.
The objects $F\in\dmcoh(X)$ are the complexes whose cohomology sheaves
$\ch^n(F)$ are coherent for all $n$, and vanish if $n\gg0$.
And $\Dqc(X)^b_c$ is also explicit: it is our old friend,
the category
traditionally denoted $\dcoh(X)$---we first met
$\dcoh(X)$ in Example~\ref{E1.100372}(iv), and it figures prominently
in the statement of Theorem~\ref{T29.1}.
\ermk

\rmk{R29.5000}
In Remark~\ref{R29.25} we observed that the category $\dcoh(X)$ has an
intrinsic description as a subcategory of $\ct=\Dqc(X)$,
it is $\ct^b_c$. The category $\dperf X$ also has an intrinsic
description, it's the subcategory $\ct^c$ of all compact
objects in $\ct=\Dqc(X)$, see Notation~\ref{N1.107589}.
With $\ct=\Dqc(X)$, 
Theorem~\ref{T29.1} is a statement about the categories
$\ct^c=\dperf X$ and $\ct^b_c=\dcoh(X)$---and it turns
out to be a special case
of the following, infinitely more general assertion.
\ermk

\thm{T29.27}
Let $R$ be a noetherian ring, and let
$\ct$ be an $R$--linear, approximable triangulated category.
Suppose there exists in $\ct$ a compact generator $G$ so that
$\Hom\big(G,G[n]\big)$ is a finite $R$--module for all $n\in\zz$.
Consider the two functors
\[
\cy:\ct^-_c\la\Hom_R^{}\big([\ct^c]\op\,,\,\Mod R\big),\qquad
\wt\cy:\big[\ct^-_c\big]\op\la\Hom_R^{}\big(\ct^b_c\,,\,\Mod R\big)
\]
defined by the formulas $\cy(B)=\Hom(-,B)$ and
$\wt\cy(A)=\Hom(A,-)$. Note that, in these formulas,
we permit all $A,B\in\ct^-_c$. But the $(-)$ in the formula
$\cy(B)=\Hom(-,B)$ is assumed to belong to $\ct^c$,
whereas the $(-)$ in the formula $\wt\cy(A)=\Hom(A,-)$ must
lie in $\ct^b_c$.
Now consider the following composites
\[\xymatrix@C+20pt@R-20pt{
\ct^b_c \ar@{^{(}->}[r]^i & \ct^-_c
\ar[r]^-{\cy} &
\Hom_R^{}\big([\ct^c]\op\,,\,\Mod R\big) \\
\big[\ct^c\big]\op \ar@{^{(}->}[r]^{\wi} & \big[\ct^-_c\big]\op
\ar[r]^-{\wt\cy} &
\Hom_R^{}\big(\ct^b_c\,,\,\Mod R\big)
}\]
We assert:
\be
\item
  The functor $\cy$ is full, and the essential image consists
  of the
  locally finite homological functors [see Explanation~\ref{X29.102345}
    for the definition of locally finite functors]. The composite $\cy\circ i$
  is fully faithful, and the essential image consists of the
  finite homological functors.
\item
  Assume there exists an integer $N>0$ and an object $G'\in\ct^b_c$
  with $\ct=\ogen {G'}_N^{(-\infty,\infty)}$.
  Then
  the functor $\wt\cy$ is full, and the essential image consists
  of the
  locally finite homological functors. The composite $\wt\cy\circ \wi$
  is fully faithful, and the essential image consists of the
  finite homological functors.
\ee
\ethm

\xpl{X29.102345}
In the statement of Theorem~\ref{T29.27}, the \emph{locally finite}
functors $\big[\ct^c\big]\op\la\Mod R$ are those functors $H$ such that
\be
\item
$H\big(A[i]\big)$ is a finite $R$--module for every $i\in\zz$ and every $A\in\ct^c$.
\item
For fixed $A\in\ct^c$ we have $H\big(A[i]\big)=0$ if $i\ll0$.
\setcounter{enumiv}{\value{enumi}}
\ee
while the locally finite functors $H:\ct^b_c\la\Mod R$ are those satisfying
\be
\setcounter{enumi}{\value{enumiv}}
\item
$H\big(A[i]\big)$ is a finite $R$--module for every $i\in\zz$ and every $A\in\ct^b_c$.
\item
For fixed $A\in\ct^b_c$ we have $H\big(A[i]\big)=0$ if $i\ll0$.
\setcounter{enumiv}{\value{enumi}}
\ee
The careful reader will observe that these definitions \emph{aren't}
dual. The finiteness of $H\big(A[i]\big)$ for every $A$ and every $i$
is, of course, obviously self-dual. But the vanishing isn't.
We might be tempted to unify the definitions into
\be
\setcounter{enumi}{\value{enumiv}}
\item
Let $\cs$ be a triangulated category. A homological functor
$H:\cs\la\Mod R$ is locally finite if $H\big(A[i]\big)$ is a finite  
$R$--module for every $i\in\zz$ and every $A\in\cs$, and for fixed
$A$ it vanishes when $i\ll0$.
\setcounter{enumiv}{\value{enumi}}
\ee
But this definition is wrong for $\big[\ct^c\big]\op$,
because its suspension functor is the negative of that of $\ct^c$.

The way to unify the two definitions is to think of them as continuity
with respect to a metric. This might become a little clearer in
the next section.
\expl

\rmk{R29.973285}
From what we have said already it's clear that the statement
about $\cy$ in
Theorem~\ref{T29.1} is a special case of Theorem~\ref{T29.27}(i).
To deduce the assertion of Theorem~\ref{T29.1}
about the functor $\wt\cy$, from Theorem~\ref{T29.27}(ii),
we need to know that the category $\dcoh(X)$ contains an object
$G'$ with $\Dqc(X)=\ogen {G'}_N^{(-\infty,\infty)}$. This is
proved in \cite[Theorem~2.3]{Neeman17}.
\ermk

\skt{S29.29}
We should say something about the proof of
Theorem~\ref{T29.27}---to keep the discussion focused let us restrict
our attention to the proof of Theorem~\ref{T29.27}(i).
For the purpose of this discussion let us fix
a compact generator $G\in\ct$ and a \tstr\
$\big(\ct^{\leq0},\ct^{\geq0}\big)$ in the
preferred equivalence class. Proposition~\ref{P29.17}
tells us that we may choose an integer $A>0$ so that
the properties of Definition~\ref{D4.1}~(i) and (ii) hold.
An easy induction on the integer $m$ leads to the following
consequence of Definition~\ref{D4.1}(ii):
\be 
\item
For every integer $m>0$ and every object $F\in\ct^{\leq0}$ there exists
a triangle $E_m\la F\la D_m\la E[1]$ with $D_m\in\ct^{\leq-m}$ and
$E\in\ogenu G{mA}{1-m-A,A}$.
\setcounter{enumiv}{\value{enumi}}
\ee 
This much is easy. Not quite so straightforward is the following:
\be 
\setcounter{enumi}{\value{enumiv}}
\item
There exists an integer $B$, depending only on $A$, with the following
property. For any integer $m>0$ and any object $F\in\ct^{\leq0}\cap\ct^-_c$
there exists a triangle $E_m\la F\la D_m\la E[1]$ with $D\in\ct^{\leq-m}$ and
$E\in\genu G{mB}{1-m-B,B}$.
\item
In fact more is true: the objects $E_m$, in either (i) or (ii) above,
can be chosen to form a sequence $E_1\la E_2\la E_3\la\cdots$ mapping to
$F$, and such that $F$ is the homotopy colimit of the sequence. It is in this
sense that the Introduction should be understood: we have
expressed $F$ as the homotopy
colimit of the (directed) Cauchy sequence $\{E_m\}$.
\ee
The reader might wish to go back to Examples~\ref{E3.7} and \ref{E3.9},
in which we explicitly worked out what the abstract theory comes
down to in the special case where $\ct=\D(R)$, the \tstr\ is the standard
one, and the compact generator $G$ is the object $R\in\D(R)$. In the
terminology of (i) and (ii) above, Examples~\ref{E3.7} and \ref{E3.9}
amount to showing that $A=B=1$ works for the special case.

Now back to the proof of Theorem~\ref{T29.27}(i). The fact that $\cy$ is
full on the category $\ct^-_c$ and fully faithful on the category $\ct^b_c$
turns out to be a straightforward consequence of (iii) above. It remains
to show that the essential image of $\cy$ is as claimed. One
containment is easy: the
fact that the essential image is contained in the locally finite (respectively
finite) functors follows directly from the the hypothesis
that $\ct$ has a compact generator $G$ so that
$\Hom\big(G,G[n]\big)$ is a finite $R$--module for all $n\in\zz$,
coupled with the fact that approximability implies the vanishing
of $\Hom\big(G,G[n]\big)$ for $n\gg0$.

It remains to prove the opposite containment. Fix a
locally finite homological functor $H:(\ct^c)\op\la\Mod R$.
We need to exhibit an object $F\in\ct^-_c$ and an isomorphism
$H\cong\cy(F)$. This actually suffices: it's easy to show that if
$\cy(F)$ is finite then $F\in\ct^b_c$.

Modifying
an old idea of Adams~\cite{Adams71} one can produce, for each
integer $m>0$, an object $F_m\in\ct^-_c$, a morphism $\cy(F_m)\la H$,
and show that this morphism is surjective when restricted to
$\cup_n\genu Gm{-n,n}\subset\ct^c$. Shifting $H$ if necessary, we
may furthermore assume that $F_m\in\ct^-_c\cap\ct^{\leq0}$ for all
$m>0$. By (ii) above we may, for each integer $m$, construct a triangle
$E_m\la F_m\la D_m\la E_m[1]$ with $E_m\in\genu G{mB}{1-m-B,B}$ and
$D_m\in\ct^{\leq-m}$. One then needs to show the existence of an
increasing sequence of integers $\{m_1^{}<m_2^{}<m_3^{}<\cdots\}$ such that
there is a Cauchy sequence $E_{m_1^{}}^{}\la E_{m_2^{}}^{}\la E_{m_3^{}}^{}\la \cdots$
converging in $\ct^-_c$ to an object $E$ with a surjection $\cy(E)\la H$.

This is the hard part. Once we have produced a surjection $\cy(E)\la H$,
one performs some tricks to deduce an isomorphism $\cy(F)\cong H$.
\eskt

\section{The categories $\dperf X$ and $\dcoh(X)$ determine each
  other}
\label{S205}

It remains to discuss Fact~\ref{F5.5}(v). We remind the reader: this is
the assertion that $\dperf X$ and $\dcoh(X)$ determine each other, as
triangulated categories.

\hst{H205.1}
Probably the first result in this direction deals with the case of
affine schemes and may be found in
Rickard~\cite[Theorem~6.4]{Rickard89b}. Rickard's result tells us that, if
$R$ and $S$ are noetherian rings, then
\[
\D^b(\proj R)\cong \D^b(\proj S) \quad\Longleftrightarrow\quad
\D^b(\mod R)\cong \D^b(\mod S)
\]
The way Rickard's theorem was proved was by analyzing the triangulated
equivalences---in other words Rickard
developed a Morita theory, and showed that
the data that produces an equivalence on the left is the same as the
data producing an equivalence on the right. Thus the proof does guarantee that
a triangulated equivalence on the left will produce a triangulated equivalence
on the right, and vice versa. But the new result is better in giving explicit
recipes for passing back and forth between $\D^b(\proj R)$ and $\D^b(\mod R)$;
there is only one scheme and two derived categories, not two schemes and four
derived categories.
Moreover: the new result doesn't assume the schemes to be affine.

A different result by Rouquier~\cite[Remark~7.50]{Rouquier08}
asserts that, if $X$ and $Y$ are projective over a field $k$, then
\[
\dcoh (X)\cong \dcoh (Y) \quad\Longleftrightarrow\quad
\dperf X\cong \dperf Y
\]
And we know both the result and the proof from the previous section:
the category $\dcoh(X)$ can be described as the category of finite homological
functors on $\big[\dperf X\big]\op$, and the category
$\big[\dperf X\big]\op$ can
be described as the finite homological functors on $\dcoh(X)$. This time
there is no Morita theory as background, we explicitly know how to
construct the categories out of each other. The drawbacks of the result,
compared to the new one, are
\be
\item
In Rouquier's statement the schemes were assumed proper over a field, and even
the improved representability theorems of the last section only work under
the assumption of properness over some commutative, noetherian ring $R$. 
\item
The triangulated structure isn't obvious. The category of (finite) homological
functors on an $R$--linear triangulated
category doesn't usually carry a triangulated
structure.
\ee
Thus the old results only work if the schemes are either affine or
proper---whereas the new result has no such restrictions---and even in
the special circumstances where one of the old results holds, the
conclusion is less powerful than what's in Fact~\ref{F5.5}(v).
\ehst

\dis{D205.3}
In the Introduction we already mentioned that a heuristic way
to think of approximability is to consider the ``metric'' defined by
a \tstr\ $\big(\ct^{\leq0},\ct^{\geq0}\big)$, and maybe study the limits
of Cauchy sequences with
respect to this metric. The reader is referred
back to Discussion~\ref{D0.1}.
Phrased in this language,  the category $\ct^-_c$
can be defined as the closure
in $\ct$ of the subcategory $\ct^c$ with respect to the
metric---see Fact~\ref{F0.9}(iv). This
suggests the recipe for constructing $\ct^b_c$
out of $\ct^c$, it remains to flesh it out a little.
\edis

\dfn{D205.5}
Let $\cs$ be a triangulated category. A \emph{metric} on $\cs$ is
a sequence of additive subcategories
$\{\cm_i,\,i\in\nn\}$, satisfying
\be
\item
$\cm_{i+1}\subset\cm_i$ for every $i\in\nn$.
\item
$\cm_i*\cm_i=\cm_i$, with the notation as in Reminder~\ref{R3.1}.
\ee
A metric $\{\cm_i\}$ is declared to be \emph{finer than the metric}
$\{\cn_i\}$
if, for every integer $i>0$, there exists an integer $j>0$ with
$\cm_j\subset\cn_i$; we denote this partial order by
$\{\cm_i\}\preceq\{\cn_i\}$.
The metrics $\{\cm_i\}$, $\{\cn_i\}$ are \emph{equivalent} if
$\{\cm_i\}\preceq\{\cn_i\}\preceq\{\cm_i\}$.
\edfn

\exm{E205.7}
Suppose $\ct$ is an approximable triangulated category, and
let $\big(\ct^{\leq0},\ct^{\geq0}\big)$ be a \tstr\ in the preferred equivalence
class.
Out of this data we can construct two examples of $\cs$'s with metrics:
\be
\item
Let $\cs$ be the subcategory $\ct^c\subset\ct$, and
put $\cm_i=\ct^c\cap\ct^{\leq-i}$.
\item
Let $\cs$ be the subcategory $\big[\ct^b_c\big]\op$, and put
$\cm\op_i=\ct_c^b\cap\ct^{\leq-i}$.
\ee
It's obvious that equivalent {\it t}--structures define equivalent metrics.
Thus up to equivalence we have a canonical metric on $\ct^c$
and a canonical metric on $\big[\ct^b_c\big]\op$. But the definition
depends on the embedding into $\ct$, which is the category
with the \tstr.
\eexm

\dfn{D205.9}
Let $\cs$ be a triangulated category with a metric $\{\cm_i\}$.
A \emph{Cauchy sequence} in $\cs$ is a sequence
$E_1\la E_2\la E_3\la\cdots$ such that, for every
pair of integers $i>0$, $j\in\zz$,
there exists
an integer $M>0$ such that,  in any triangle 
$E_m\la E_{m'}\la D_{m,m'}$ with $M\leq m<m'$, the object $\T^jD_{m,m'}$ lies
in $\cm_i$.  
\edfn

\nin
It is clear that the Cauchy sequences depend only on the equivalence
class of the metric.

\con{C205.11}
Now suppose $\cs$ is an essentially small triangulated category with a metric.
If we write $\MMod\cs$ for the category of additive functors
$\cs\op\la\ab$, then the Yoneda functor is a fully faithful embedding
$Y:\cs\la\MMod\cs$. We remind the reader: the formula is $Y(A)=\Hom(-,A)$.

In the abelian
category $\MMod\cs$ we can form colimits. We define
$\fl(\cs)\subset\MMod\cs$
to be the full subcategory of all colimits in $\MMod\cs$ of
Cauchy sequences in $\cs$, and define $\fs(\cs)\subset\fl(\cs)$
to be the full subcategory
\[
\fs(\cs)\eq\fl(\cs)\bigcap\left[\bigcap_{j\in\zz}\,\,\bigcup_{i\in\nn}\big[Y(\T^j\cm_i)\big]^\perp\right]
\]
Recall: an object $M\in\MMod\cs$ belongs to
$\big[Y(\T^j\cm_i)\big]^\perp$ if,
for every object $m\in\T^j\cm_i$, we have $\Hom\big(Y(m),M\big)=0$.

It's clear that the subcategories $\fs(\cs)\subset\fl(\cs)\subset\MMod\cs$
depend only on the equivalence class of the metric.
\econ

The first result, which may be found in \cite[Theorem~2.11]{Neeman18A},
asserts:

\thm{T205.13}
The category $\fs(\cs)$ is a triangulated category, with the triangles
being the colimits in $\MMod\cs$ of Cauchy sequences of triangles in $\cs$.
\ethm

\nin
Thus the triangulated structure on $\fs(\cs)$ also depends only on the
equivalence class of the metric.

And the second result, which is in
\cite[Example~4.2 and Proposition~5.6]{Neeman18A}, gives:

\thm{T205.15}
With the metrics as in Example~\ref{E205.7}~(i) and (ii), we have
triangulated equivalences
\be
\item
$\fs(\ct^c)=\ct^b_c$.
\item
If $\ct$ is noetherian then
$\fs\big(\big[\ct^b_c\big]\op\big)=\big[\ct^c\big]\op$.
\ee
\ethm

\ntn{N205.17}
We owe the reader an explanation of the hypothesis in
Theorem~\ref{T205.15}(ii). Let $\ct$ be an approximable
triangulated category, and
let $\big(\ct^{\leq0},\ct^{\geq0}\big)$ be a \tstr\ in the preferred equivalence
class. The category $\ct$ is \emph{noetherian} if there exists an
integer $N>0$ so that, for every object $F\in\ct^-_c$, there exists
a triangle $F'\la F\la F''$ in $\ct^-_c$ with
$F'\in\ct^-_c\cap\ct^{\leq N}$ and $F''\in\ct^-_c\cap\ct^{\geq0}$.

It's clear that replacing the \tstr\ $\big(\ct^{\leq0},\ct^{\geq0}\big)$ by
an equivalent one will only have the effect of changing the integer $N$;
the definition doesn't depend on the choice of \tstr\ in the
preferred equivalence class.
\entn

\rmk{R205.19}
The way the metrics were presented, in Example~\ref{E205.7}, depends
on the embedding of $\ct^c$ and $\ct^b_c$ in $\ct$. After all they are
defined in terms of the preferred equivalence class of {\it t}--structures,
which makes sense only in $\ct$.

But there are recipes that cook up equivalent metrics directly from
$\ct^c$ and $\ct^b_c$. The reader is referred to
\cite[Remark~4.7 as well as Propositions~4.8 and 6.5]{Neeman18A}.
\ermk

\rmk{R205.21}
In this survey we have often said that, so far, the main applications
of the theory have been to algebraic geometry---it's the fact
that $\Dqc(X)$ is approximable which has had the far-reaching consequences
to date.
In this context Theorem~\ref{T205.15} is the first exception, it has
striking consequences for categories not coming from
algebraic geometry. From
the theorem we learn the following (among other things).
 \be
    \item
      Let $R$ be any ring, possibly noncommutative. The recipe takes the triangulated category $\D^b(\proj R)$ and out of it constructs the triangulated category $\D^-(\proj R)\cap\D^b(\Mod R)$. The objects of this intersection
      are the bounded-above cochain complexes of finitely generated, projective
      $R$--modules, with bounded cohomology.

      If $R$ is a coherent ring this category is equivalent to $\D^b(\mod R)$.
    \item If $R$ is a coherent (possibly noncommutative) ring,
      then the recipe takes $\big[\D^b(\mod R)\big]\op$ and out of it
      constructs $\big[\D^b(\proj R)\big]\op$.
    \item Out of the homotopy category of finite spectra we construct the homotopy category of spectra with finitely many nonzero stable homotopy groups, all of them finitely generated.
      \item 
      Out of the homotopy category of spectra with finitely many nonzero stable homotopy groups, all of them finitely generated, we construct the homotopy category of finite spectra.
  \ee
\ermk

\rmk{R205.91763}
When I posted the current article on the archive Steve Lack wrote to inform me
that metrics in categories aren't new. From Richard Garner I later learned
that completing with respect to such metrics has also been done before. Hence
I wrote the survey~\cite{Neeman19}, which offers a much-expanded treatment
of the material discussed in Section~\ref{S205} and includes some remarks 
about the category-theory literature that preceded.
\ermk

\appendix

\section{Some dumb maps in $\D_\fc^{\fc'}(\ca)$, and the proof that
  the third map of the triangle is a cochain map} 
\label{A1}

For any object $A$ in the abelian category $\ca$,
we will write $\wt A$ for the cochain
complex
\[\xymatrix@C+20pt{
\cdots \ar[r]& 0  \ar[r]& A \ar@{=}[r]& A  \ar[r]
&  0 \ar[r] &0 \ar[r] &\cdots
}\]
with $A$ in degrees $-1$ and 0. Since the cohomology of the complex $\wt A$
vanishes, in the derived category $\D_\fc^{\fc'}(\ca)$
the morphism $0\la \wt A$ is
an isomorphism. We reiterate: $\wt A$ is
nothing more than a complicated representative of
the isomorphism class of the zero object in $\D_\fc^{\fc'}(\ca)$.
With the conventions of Example~\ref{E1.7} and Notation~\ref{N1.8},
for any integer $i\in\zz$ we may form
the object $\wt A[i]$, which is the cochain complex
\[\xymatrix@C+20pt{
\cdots \ar[r]& 0  \ar[r]& A \ar[r]^-{\e}& A  \ar[r]
&  0 \ar[r] &0 \ar[r] &\cdots
}\]
with $A$ in degrees $-i-1$ and $-i$, and where $\e=(-1)^i$.

Let $X^*\in\D_\fc^{\fc'}(\ca)$ be any object, meaning $X^*$ is a cochain complex 
\[\xymatrix@C+17pt{
\cdots \ar[r]& X^{-2}  \ar[r]^-{\partial_X^{-2}}& X^{-1}  \ar[r]^-{\partial_X^{-1}}& X^{0}  \ar[r]^-{\partial_X^{0}}
&  X^{1}\ar[r]^-{\partial_X^{1}}& X^{2} \ar[r] &\cdots
}\]
The cochain maps $X^*\la\wt A[-i]$ are in bijection with morphisms
$\theta^i:X^{i}\la A$ in $\ca$, the bijection takes $\theta^i$ to the cochain map
\[\xymatrix@C+15pt{
\cdots \ar[r]& X^{i-2}  \ar[d]\ar[r]^-{\partial_X^{i-2}}& X^{i-1} \ar[d]_-{\e\theta^i\partial_X^{i-1}} \ar[r]^-{\partial_X^{i-1}}& X^{i}\ar[d]^{\theta^i}  \ar[r]^-{\partial_X^{i}}
&  X^{i+1} \ar[d]\ar[r]^-{\partial_X^{i+1}}& X^{i+2}\ar[d] \ar[r] &\cdots\\
\cdots \ar[r]& 0  \ar[r]& A  \ar[r]_{\e}& A  \ar[r]
&  0 \ar[r] &0 \ar[r] &\cdots
}\]
Given two cochain complexes $X^*$ and $Y^*$,
as well as a morphism $\theta^i:X^i\la Y^{i-1}$ in $\ca$,
we may form the composite
\[\xymatrix@C+15pt{
\cdots \ar[r]& X^{i-2}  \ar[d]\ar[r]^-{\partial_X^{i-2}}& X^{i-1} \ar[d]_-{\e\partial_X^{i-1}} \ar[r]^-{\partial_X^{i-1}}& X^{i}\ar[d]^{\id}  \ar[r]^-{\partial_X^{i}}
&  X^{i+1} \ar[d]\ar[r]^-{\partial_X^{i+1}}& X^{i+2}\ar[d] \ar[r] &\cdots\\
\cdots \ar[r]& 0  \ar[r]\ar[d]& X^{i}  \ar[r]_{\e}\ar[d]_{\e\theta^i}& X^{i}  \ar[r]
\ar[d]^-{\partial_Y^{i-1}\theta^i}
&  0 \ar[r]\ar[d] &0 \ar[r]\ar[d] &\cdots \\
\cdots \ar[r]& Y^{i-2}  \ar[r]^-{\partial_Y^{i-2}}& Y^{i-1}  \ar[r]^-{\partial_Y^{i-1}}& Y^{i} \ar[r]^-{\partial_Y^{i}}
&  Y^{i+1}\ar[r]^-{\partial_Y^{i+1}}& Y^{i+2} \ar[r] &\cdots
}\]
which is manifestly a cochain map; we will
denote it by $\wt\theta^i$. Of course in the category $\D_\fc^{\fc'}(\ca)$ the
cochain map $\wt\theta^i$ must
be equal to the zero map, all we have done is produced many representatives
of the zero map. We may also sum over $i\in\zz$; given any collection
of morphisms $\theta^i:X^i\la Y^{i-1}$ in $\ca$ we may form
$\oplus_{i=-\infty}^\infty\wt\theta^i$, which is a cochain map
$\wt\theta^*:X^*\la Y^*$. Needless to say the map $\wt\theta^*$ also vanishes
in $\D_\fc^{\fc'}(\ca)$.

Now go back to Example~\ref{E1.7}: in the example we construct a diagram
\[\xymatrix@C+8pt{
\cdots \ar[r]& Z^{-2}  \ar[d]^-{h^{-2}}\ar[r]^-{\partial_Z^{-2}}& Z^{-1} \ar[d]^-{h^{-1}} \ar[r]^-{\partial_Z^{-1}}& Z^{0}\ar[d]^-{h^{0}}  \ar[r]^-{\partial_Z^{-0}}
&  Z^{1} \ar[d]^-{h^{1}}\ar[r]^-{\partial_Z^{1}}& Z^{2}\ar[d]^--{h^{2}} \ar[r] &\cdots\\
\cdots \ar[r]& X^{-1}\ar[d]^-{f^{-1}}  \ar[r]_-{-\partial_X^{-1}}& X^{0}
\ar[d]^-{f^{0}}\ar[r]_-{-\partial_X^{0}}& X^{1}
\ar[d]^-{f^{1}} \ar[r]_-{-\partial_X^{1}}
&  X^{2}\ar[d]^-{f^{2}}  \ar[r]_-{-\partial_X^{2}}& X^{3}
\ar[d]^-{f^{3}} \ar[r] &\cdots\\
\cdots \ar[r]& Y^{-1}  \ar[r]_-{-\partial_Y^{-1}}& Y^{0}  \ar[r]_-{-\partial_Y^{0}}& Y^{1}  \ar[r]_-{-\partial_Y^{1}}
&  Y^{2} \ar[r]_-{-\partial_Y^{2}}& Y^{3} \ar[r] &\cdots
}\]
The construction is such that, if we delete the middle row, then the composite
is $\wt\theta^*$ for the explicit $\theta^*$ chosen in the Example. In
particular: deleting the middle row gives a commutative diagram. Deleting
the top row yields a commutative diagram, we are left with nothing other
than the given cochain map $f^*[1]$. We conclude that, for each $i$,
the composites
from top left to bottom right in the diagram
\[\xymatrix@C+8pt{
Z^{i}\ar[d]^-{h^{i}}  \ar[r]_-{\partial_Z^{i}}
&  Z^{i+1} \ar[d]^-{h^{i+1}}\\
X^{i+1}
\ar[r]_-{-\partial_X^{i+1}}
&  X^{i+2}\ar[d]^-{f^{i+2}} \\
&  Y^{i+2} 
}\]
must agree. Since the map $f^{i+2}$ is a (split) monomorphism the square
commutes, and we conclude that the diagram
\[\xymatrix@C+8pt{
\cdots \ar[r]& Z^{-2}  \ar[d]^-{h^{-2}}\ar[r]^-{\partial_Z^{-2}}& Z^{-1} \ar[d]^-{h^{-1}} \ar[r]^-{\partial_Z^{-1}}& Z^{0}\ar[d]^-{h^{0}}  \ar[r]^-{\partial_Z^{-0}}
&  Z^{1} \ar[d]^-{h^{1}}\ar[r]^-{\partial_Z^{1}}& Z^{2}\ar[d]^--{h^{2}} \ar[r] &\cdots\\
\cdots \ar[r]& X^{-1}  \ar[r]_-{-\partial_X^{-1}}& X^{0}  \ar[r]_-{-\partial_X^{0}}& X^{1}  \ar[r]_-{-\partial_X^{1}}
&  X^{2} \ar[r]_-{-\partial_X^{2}}& X^{3} \ar[r] &\cdots
}\]
commutes. Thus $h^*$ is indeed a cochain map, as promised in Example~\ref{E1.7}.

\section{The assumption that the short exact sequences of cochain complexes
are degreewise split is harmless}
\label{A2}

Given a cochain complex $X^*$,
an object $A\in\ca$ and
a morphism $\theta^i:X^i\la A$,
in Appendix~\ref{A1} we constructed a corresponding cochain map
$X^*\la\wt A[-i]$. In the special case where $A=X^i$ and $\theta^i$ is the
identity $\id:X^i\la X^i$ the general recipe specializes
to the cochain map $\rho^i$ below
\[\xymatrix@C+15pt{
\cdots \ar[r]& X^{i-2}  \ar[d]\ar[r]^-{\partial_X^{i-2}}& X^{i-1} \ar[d]_-{\e\partial_X^{i-1}} \ar[r]^-{\partial_X^{i-1}}& X^{i}\ar@{=}[d]  \ar[r]^-{\partial_X^{i}}
&  X^{i+1} \ar[d]\ar[r]^-{\partial_X^{i+1}}& X^{i+2}\ar[d] \ar[r] &\cdots\\
\cdots \ar[r]& 0  \ar[r]& X^i  \ar[r]_{\e}& X^i  \ar[r]
&  0 \ar[r] &0 \ar[r] &\cdots
}\]
For this section the key point is that, in degree $i$, the cochain map
$\rho^i:X^*\la \wt X^i[-i]$ is a split
monomorphism. Taking the direct sum over 
$i$
produces
\[\xymatrix@C+40pt{
X^*\ar[r] & \displaystyle\bigoplus_{i\in\zz}\wt X^i[-i]
}\]
We denote this map $\rho^*:X^*\la\wt X^*$ and observe that
\be
\item
The object $\wt X^*=\oplus_{i\in\zz}\wt X^i[-i]$ vanishes in $\D_\fc^{\fc'}(\ca)$.
\item
The cochain map $\rho^*$ is a split monomorphism in each degree.
\ee
Now suppose we're given a
short exact sequence of cochain complexes
$X^*\stackrel{f^*}\la Y^*\stackrel{g^*}\la Z^*$. We may form
the commutative diagram of cochain complexes, where
the rows are exact
\[\xymatrix@C+10pt{
0\ar[r] &X^*\ar[rr]^-{\left(\begin{array}{c}
\rho^* \\ f^* 
\end{array}\right)}
\ar@{=}[d]& & \wt X^*\oplus Y^* \ar[d]^-{\pi}
\ar[rr]^-{\wt g^*} & & \wt Z^*\ar[d]^-{\ph^*}\ar[r] & 0\\
0\ar[r] &X^*\ar[rr]^-{f^*}&&  Y^* \ar[rr]^-{g^*} && Z^* \ar[r]& 0
}\]
We know that the vertical
maps $\id:X\la X$ and $\pi:\wt X^*\oplus Y^*\la Y^*$
both induce isomorphisms in cohomology. The 5-lemma,
applied to the long exact sequences in cohomology
that come
from the short exact sequences of cochain complexes
in the rows, tells us that
the vertical morphism $\ph^*$ also induces
an isomorphism in cohomology. Hence the top row is isomorphic in
$\D_\fc^{\fc'}(\ca)$ to the bottom row, and the top row is degreewise split.
This is the sense in which we said, back in Example~\ref{E1.7}, that
up to isomorphism in $\D_\fc^{\fc'}(\ca)$ we may assume our short exact sequence
of cochain complexes is
degreewise split.

\section{Translating the approach to derived categories
  we presented here to the more
  standard one in the literature}
\label{A3}

Our presentation of derived categories has
been minimalist---we have tried to be accurate
without providing any information that wasn't absolutely
indispensable. This means
that the student who is seeing this for the first time, and would like
to look for more detail elsewhere in the literature, might have a hard time
reconciling what's here with other, more expansive accounts. This appendix
was written to help.

\dfn{DA3.-3}
Suppose $\ca$ is an abelian category. Following tradition
we define the two categories
$\CC^{\fc'}_\fc(\ca)$ and $\D^{\fc'}_\fc(\ca)$
to have the same objects: cochain complexes in $\ca$ subject to
identical restrictions. The morphisms are
\be
\item
In the category $\CC^{\fc'}_\fc(\ca)$ a morphism is a cochain map.
\item
We have already met $\D^{\fc'}_\fc(\ca)$: it is obtained from
$\CC^{\fc'}_\fc(\ca)$ by formally inverting the morphisms in 
$\CC^{\fc'}_\fc(\ca)$ inducing isomorphisms in cohomology.
\ee
\edfn

\rmk{RA3.1}
We give the following list, translating the constructions of the previous
two appendices to more standard language.
\be
\item
Given a
pair of objects $X^*$ and $Y^*$ in $\CC^{\fc'}_\fc(\ca)$,
that is a pair of cochain complexes $X^*$ and $Y^*$,
as well as a sequence of
morphisms
$\{\theta^i:X^i\la Y^{i-1},\,i\in\zz\}$ in the category $\ca$,
Appendix~\ref{A1} constructed for us (among other things)
a morphism $\wt\theta^*:X^*\la Y^*$ in the category
$\CC^{\fc'}_\fc(\ca)$. It is traditional to say that
the cochain map $\wt\theta^*$ is \emph{null homotopic,} and the
sequence of maps $\{\theta^i:X^i\la Y^{i-1},\,i\in\zz\}$ is called a
\emph{homotopy of $\wt\theta^*$ with the zero map.} More generally: two
cochain maps $f^*,g^*:X^*\la Y^*$ are declared to be \emph{homotopic to each
  other} if there exists a homotopy of $f^*-g^*$ with the zero map. That is:
if $f^*-g^*=\wt\theta^*$ with $\wt\theta^*$ as in Appendix~\ref{A1}.
\item
In Appendix~\ref{A2} we constructed, for every cochain complex
$X^*$, another cochain complex $\wt X^*$ and a map $\rho:X^*\la \wt X^*$
in the category $\CC^{\fc'}_\fc(\ca)$.
Note that, by the construction of the map $\wt\theta^*:X^*\la Y^*$ of
(i) above (see Appendix~\ref{A1}), the map $\wt\theta^*$ comes
with a factorization $X^*\stackrel\rho\la \wt X^*\stackrel\theta\la Y^*$
in the category $\CC^{\fc'}_\fc(\ca)$.
\item
Given a cochain map $f^*:X^*\la Y^*$, the object $\wt X^*\oplus Y^*$
is isomorphic in $\CC^{\fc'}_\fc(\ca)$ to a complex
traditionally
called the \emph{mapping cylinder} of $f^*$. This is obviously dumb terminology
since $\wt X^*\oplus Y^*$ manifestly
doesn't depend on $f^*$. But the isomorph in $\CC^{\fc'}_\fc(\ca)$,
that is traditionally presented as the
mapping cylinder, seems to depend on $f^*$---needless to
to say, the isomorphism with $\wt X^*\oplus Y^*$
involves $f^*$.
\item
In Appendix~\ref{A2} we considered the degreewise split short
exact sequence of cochain complexes  
\[\xymatrix@C+10pt{
0\ar[r] &X^*\ar[rr]^-{\left(\begin{array}{c}
\rho^* \\ f^* 
\end{array}\right)}
& & \wt X^*\oplus Y^* 
\ar[rr]^-{\wt g^*} & & \wt Z^*\ar[r] & 0
}\]
The cochain complex $\wt Z^*$ is isomorphic in $\CC^{\fc'}_\fc(\ca)$
to a complex
traditionally known as the \emph{mapping cone}
on $f^*$.
\item
When $f^*:X^*\la Y^*$ is the identity $\id:X^*\la X^*$, the mapping cone
is isomorphic in $\CC^{\fc'}_\fc(\ca)$
to $\wt X^*$. Thus the object we call $\wt X^*$ is traditionally  
called the \emph{mapping cone on the identity.} Once again it is traditional
to give an isomorphic complex, which is more complicated-looking
than the $\wt X^*$ of Appendix~\ref{A2}.
\ee
\ermk

\dfn{DA3.3}
Suppose $\ca$ is an abelian category. Following tradition
we define, in addition to the two categories of Definition~\ref{DA3.-3},
yet another category
$\K^{\fc'}_\fc(\ca)$.
It has the same objects as $\CC^{\fc'}_\fc(\ca)$ or $\D^{\fc'}_\fc(\ca)$.
However:
\be
\item
In the category $\K^{\fc'}_\fc(\ca)$ a morphism is a homotopy
equivalence class of cochain maps. Thus two cochain maps are declared
to be equal in $\K^{\fc'}_\fc(\ca)$ if they're homotopic.
\ee
\edfn

\lem{LA3.3}
The natural functor $F:\CC^{\fc'}_\fc(\ca)\la\D^{\fc'}_\fc(\ca)$,
that is the universal functor out of $\CC^{\fc'}_\fc(\ca)$ which sends
every cohomology isomorphism to an isomorphism, factors uniquely through
$\K^{\fc'}_\fc(\ca)$.
\elem

\prf
Let $X^*$ be any object in $\CC^{\fc'}_\fc(\ca)$, and consider the following
morphisms in $\CC^{\fc'}_\fc(\ca)$
\[\xymatrix@C+10pt{
X^*\ar@<0.5ex>[rr]^-{\left(\begin{array}{c}
\rho^* \\ \id 
\end{array}\right)}\ar@<-0.5ex>[rr]_-{\left(\begin{array}{c}
0 \\ \id 
\end{array}\right)}
& & \wt X^*\oplus X^* 
\ar[rr]^-{(0,\id)} & & X^*
}\]
The two composites are equal, in each case the composite is the identity
map
$\id:X^*\la X^*$. But the map $(0,\id):\wt X^*\oplus X^*\la X^*$ is an
isomorphism in cohomology [because $\wt X^*$ is acyclic], therefore
$F:\CC^{\fc'}_\fc(\ca)\la\D^{\fc'}_\fc(\ca)$ takes $(0,\id)$ to an
isomorphism in $\D^{\fc'}_\fc(\ca)$. Hence the functor
$F:\CC^{\fc'}_\fc(\ca)\la\D^{\fc'}_\fc(\ca)$ must take the two cochain
maps\footnote{Recall Remark~\ref{RA3.1}(iii): for any
  cochain map $f^*:X^*\la X^*$ there is
  a cochain complex
  traditionally called the \emph{mapping cylinder of $f^*$,} and all
  it is, as an object of the
  category $\CC^{\fc'}_\fc(\ca)$, is a complicated isomorph of 
  $\wt X^*\oplus X^*$. If
  $f^*:X^*\la X^*$ is the identity map then the two morphisms
  $\xymatrix{X^*\ar@<0.3ex>[r]\ar@<-0.3ex>[r] & \wt X^*\oplus X^*}$,
  studied in the proof of Lemma~\ref{LA3.3},
  traditionally go by the name ``the inclusions of the front and back faces''
  of the mapping cylinder.

  Note that, although the mapping cylinder of $f^*$ is independent of $f^*$
  up to isomorphism in $\CC^{\fc'}_\fc(\ca)$,
  the inclusion of the back face depends on $f^*$.
}
\[\xymatrix@C+30pt{
X^*\ar@<0.5ex>[rr]^-{\left(\begin{array}{c}
\rho^* \\ \id 
\end{array}\right)}\ar@<-0.5ex>[rr]_-{\left(\begin{array}{c}
0 \\ \id 
\end{array}\right)}
& & \wt X^*\oplus X^* 
}\]
to equal morphisms in $\D^{\fc'}_\fc(\ca)$.  But then, for any homotopy
$\{\theta^i:X^i\la Y^{i-1},\,i\in\zz\}$ as in Remark~\ref{RA3.1}(i), the
two composites
\[\xymatrix@C+10pt{
X^*\ar@<0.5ex>[rr]^-{\left(\begin{array}{c}
\rho^* \\ \id 
\end{array}\right)}\ar@<-0.5ex>[rr]_-{\left(\begin{array}{c}
0 \\ \id 
\end{array}\right)}
& & \wt X^*\oplus X^* 
\ar[rr]^-{(\theta,g^*)} & & Y^*
}\]
must map under $F:\CC^{\fc'}_\fc(\ca)\la\D^{\fc'}_\fc(\ca)$ to
equal morphisms in $\D^{\fc'}_\fc(\ca)$. If $\{\theta^i:X^i\la Y^{i-1},\,i\in\zz\}$
is a homotopy with $f^*-g^*=\wt\theta^*=\theta\rho^*$ we discover that,
in the category
$\D^{\fc'}_\fc(\ca)$, we must have $F(f^*)=F(g^*)$.
\eprf

\rmk{RA3.5}
It now follows easily that $\D^{\fc'}_\fc(\ca)$ can also be obtained
from $\K^{\fc'}_\fc(\ca)$ by formally inverting the maps inducing
isomorphisms in cohomology. This brings us to the traditional approach:
one proves first that the category $\K^{\fc'}_\fc(\ca)$ is triangulated,
then formally inverts. And there is a general theorem of Verdier,
giving conditions under which the process of formally inverting
morphisms takes one triangulated category to another.
\ermk

\rmk{RA3.7}
In Example~\ref{E1.7} we started in $\CC^{\fc'}_\fc(\ca)$
with a degreewise split
short exact sequence of cochain complexes
$X^*\stackrel{f^*}\la Y^*\stackrel{g^*}\la Z^*$ and, after
choosing for every $i\in\zz$ a degreewise splitting $\theta^i:Z^i\la Y^i$,
we extended in $\CC^{\fc'}_\fc(\ca)$ to form the sequence of morphisms
$X^*\stackrel{f^*}\la Y^*\stackrel{g^*}\la Z^*\stackrel{h^*}\la \T X^*$.
The triangles in $\K^{\fc'}_\fc(\ca)$ can be defined to be
the isomorphs
in $\K^{\fc'}_\fc(\ca)$ of the 
$X^*\stackrel{f^*}\la Y^*\stackrel{g^*}\la Z^*\stackrel{h^*}\la \T X^*$
that come from the construction of Example~\ref{E1.7}.

We leave it to the reader to check (if she wishes) that the triangles
we construct are isomorphic in $\K^{\fc'}_\fc(\ca)$ to the
standard ones in the literature. The interested reader can
amuse herself by furthermore
checking that the axioms of triangulated categories
are satisfied, using either the description of the triangles in
the paragraph above or the more
standard one in the literature.

Even if the reader is feeling lazy today, the honest truth is
that the existing literature won't help---it tends to leave this
verification to the reader,
there would be no benefit in working
with the standard description of triangles. Doing this will not reduce
the amount of labor the indolent reader is asked to perform, the checking
will still be her task, only the details that need to be proved will
shift a little.

The good news is that,
once the reader has verified that the category $\K^{\fc'}_\fc(\ca)$ is
triangulated, the passage from $\K^{\fc'}_\fc(\ca)$ to
$\D^{\fc'}_\fc(\ca)$ becomes well-documented. The
theorem of Verdier, alluded to in Remark~\ref{RA3.5}, is explicit
enough to provide helpful information
about the calculus of fractions involved---and
this may be found in any of the standard treatments in the literature.
\ermk

\def\cprime{$'$}
\providecommand{\bysame}{\leavevmode\hbox to3em{\hrulefill}\thinspace}
\providecommand{\MR}{\relax\ifhmode\unskip\space\fi MR }
\providecommand{\MRhref}[2]{%
  \href{http://www.ams.org/mathscinet-getitem?mr=#1}{#2}
}
\providecommand{\href}[2]{#2}

\end{document}